\documentclass[12pt]{amsart}
\usepackage{amssymb,amsmath}
\makeatletter
\def\@cite#1#2{{\m@th\upshape\bfseries%
[{#1\if@tempswa{\m@th\upshape\mdseries, #2}\fi}]}}
\makeatother
\theoremstyle{plain}
\newtheorem{thm}{Theorem}[section]
\newtheorem{lem}[thm]{Lemma}
\newtheorem{cor}[thm]{Corollary}
\newtheorem{prop}[thm]{Proposition}
\theoremstyle{definition}
\newtheorem{rem}[thm]{Remark}

\newtheorem{defn}[thm]{Definition}

\newtheorem{eg}[thm]{Example}

\newcommand{\Prf}{\noindent\textbf{Proof.\ }}
\newcommand{\bx}{\strut\hfill$\blacksquare$\medbreak}

\newcommand{\ca}{\mathrm{C}^*}

\newcommand{\ol}{\overline}
\newenvironment{sbmatrix}{\left[\begin{smallmatrix}}{\end{smallmatrix}\right]}

\newcommand{\td}{\widetilde}

\DeclareMathOperator*{\sotlim}{\textsc{sot}--lim}

\newcommand{\wot}{\textsc{wot}}

\newcommand{\bbB}{{\mathbb{B}}}
\newcommand{\bbC}{{\mathbb{C}}}
\newcommand{\bbD}{{\mathbb{D}}}
\newcommand{\bbF}{{\mathbb{F}}}

\newcommand{\bbZ}{{\mathbb{Z}}}
 \newcommand{\A}{{\mathcal{A}}}
 \newcommand{\B}{{\mathcal{B}}}
 \newcommand{\C}{{\mathcal{C}}}
 \newcommand{\D}{{\mathcal{D}}}

\renewcommand{\H}{{\mathcal{H}}}
 \newcommand{\I}{{\mathcal{I}}}
 \newcommand{\J}{{\mathcal{J}}}
 \newcommand{\K}{{\mathcal{K}}}
\renewcommand{\L}{{\mathcal{L}}}
 \newcommand{\M}{{\mathcal{M}}}
 \newcommand{\N}{{\mathcal{N}}}
\renewcommand{\O}{{\mathcal{O}}}

\renewcommand{\S}{{\mathcal{S}}}
 \newcommand{\T}{{\mathcal{T}}}
 \newcommand{\U}{{\mathcal{U}}}
 
 \newcommand{\W}{{\mathcal{W}}}

\newcommand{\upchi}{{\raise.35ex\hbox{$\chi$}}}
\newcommand{\fA}{{\mathfrak{A}}}

\newcommand{\fL}{{\mathfrak{L}}}

\newcommand{\fR}{{\mathfrak{R}}}

\newcommand{\qand}{\quad\text{and}\quad}

\newcommand{\qfor}{\quad\text{for}\quad}
\newcommand{\qwhere}{\quad\text{where}\quad}

\newcommand{\rad}{\operatorname{rad}}
\newcommand{\badedges}{\operatorname{B}}
\newcommand{\edges}{\operatorname{E}}

\newcommand{\Ad}{\operatorname{Ad}}

\newcommand{\Alg}{\operatorname{Alg}}

\newcommand{\dist}{\operatorname{dist}}

\newcommand{\Lat}{\operatorname{Lat}}

\newcommand{\ran}{\operatorname{Ran}}
\newcommand{\rank}{\operatorname{rank}}

\newcommand{\spn}{\operatorname{span}}

\newcommand{\proj}{\operatorname{Proj}}
\newcommand{\pisom}{\operatorname{Pisom}}
\newcommand{\fngee}{\bbF^+\!(G)}

\newcommand{\bofh}{\B(\H)}

\newcommand{\mulam}{\nu_{\lambda,i}}
\newcommand{\nulam}{\nu_{\lambda,i}}

\newcommand{\flgee}{\fL_G}
\newcommand{\frgee}{\fR_G}

\begin{document}

%%%%%%%%%%%%%%%%%%%%%%%%%%%%%%%%%%%%%%%%%%%%%%
%%%%%%%%%%
\title[Free Semigroupoid Algebras]%
{Free Semigroupoid Algebras}
%\thanks{draft May, 2003}
%
\author[D.W.Kribs and S.C.Power]{David~W.~Kribs and
Stephen~C.~Power}
\thanks{{\it key words and phrases. } Hilbert space, Fock space, partial
isometry, directed graph,  free semigroupoid, commutant,
reflexive algebra,
automorphism, amalgamated free product, partly free algebra.}
\thanks{first author partially supported by a Canadian NSERC
Post-doctoral Fellowship.}
\address{Department of Mathematics and Statistics, University of Guelph,
Guelph, Ontario, Canada  N1G 2W1} \email{dkribs@uoguelph.ca}
\address{Department of Mathematics and Statistics, Lancaster University, Lancaster, England
LA1 4YW}
\email{s.power@lancaster.ac.uk}

\date{}
\begin{abstract}
Every countable directed graph generates a Fock space Hilbert
space and a family of partial isometries.
These operators also arise from the left regular representations of
free semigroupoids derived from directed graphs. We develop a structure theory for the weak
operator topology closed algebras generated by these representations,
which we call {\it free semigroupoid algebras}.
We characterize semisimplicity in terms of the graph and show
explicitly in the case of finite graphs how the Jacobson radical is determined.
We provide a diverse collection of examples including; algebras with
free behaviour, and examples which can be represented as matrix function
algebras. We show how these algebras can be presented and decomposed in terms of
amalgamated free products.
We determine the commutant, consider invariant subspaces, obtain a Beurling
theorem for them, conduct an eigenvalue
analysis, give an
elementary proof of reflexivity, and discuss hyper-reflexivity. Our main theorem
shows the graph to be a complete unitary
invariant for the algebra. This classification theorem makes use of
an analysis of unitarily implemented automorphisms.
We give a graph-theoretic description of when these algebras are {\it partly
free}, in the sense that they contain a copy of a free semigroup algebra.
\end{abstract}
\maketitle
%%%%%%%%%%%%%%%%%%%%%%%%%%%%%%%%%%%%%%%%%%%%%%
%%%%%%%%%%%

%\setcounter{tocdepth}{1}
%\tableofcontents

%%%%%%%%%%%%%%%%%%%%%%%%%%%%%%%%%%%%%%%%%%%%%%
%%
\section{Introduction}\label{S:realintro}
%%%%%%%%%%%%%%%%%%%%%%%%%%%%%%%%%%%%%%%%%%%%%%
%%

We initiate the study of a new class of operator algebras which we
call {\it free semigroupoid algebras}. These algebras include, as
special cases, the space $H^\infty$ realized as the  {\it analytic
Toeplitz algebra} \cite{Douglas_text,Hoffman_text} and the
prototypical {\it free semigroup algebras} $\fL_n$ of Popescu,
Davidson and Pitts (see
\cite{AP,DKP,DP1,DP2,Kribs_factor,Pop_fact,Pop_beur}). However our
context also embraces finite-dimensional operator algebras
(inflation algebras), finite and infinite matrix function
algebras, as well as operator algebras with free or partly free
structure. Thus we have a unifying framework for the free
semigroup algebras which includes the classical analytic Toeplitz
algebra as part of a general class, rather than as an exceptional
case. The framework also includes various nest algebras and finite
dimensional digraph algebras but in a represented form wherein the
commutant algebras have a similar character.

The generators of these algebras are families of partial
isometries and projections which
arise from countable directed graphs. Each graph $G$ naturally determines a
generalized Fock space Hilbert space and partial creation operators which
act
on the space. Alternatively, these operators come from the left regular
representation of the {\it free semigroupoid} derived from the
directed graph,
thus giving further credence to the terminology.
%Our principal focus will be on the algebras $\flgee$ such that for each partial isometry
%generator
%$S$ the
%following one-dimensional defect condition holds:
%\[
%\dim \Big( S^*S - \sum_{S_jS_j^* \leq S^*S} S_jS_j^* \Big) = 1.
%\]

The norm closed versions of these  algebras, in the case of finite graphs, were
considered by Muhly \cite{Muhly_fd} with particular reference to commutant lifting, the Shilov
boundary and the $\ca$-envelope.
Subsequent work with Solel \cite{MS_CJM} on more general (norm-closed)
tensor algebras addressed the structure of ideals, Wold
decompositions and a Beurling type invariant subspace theorem. There is some
overlap here with our development which we point out in Sections~4
and 8.

Presently there is considerable interest in $\ca$-algebras generated by families
of partial isometries associated with directed graphs and how these algebras
relate to the Elliott classification programme. The set of generators in this case is sometimes
referred to as a
Cuntz-Krieger $E$-family (for instance see  \cite{BHRS,Kumjian1,Kumjian2}).
On the other hand, our generators are of Cuntz-Krieger-Toeplitz type and we shall
see that it is natural to include the initial projections. Moreover, our analysis
is thoroughgoingly non-selfadjoint and spatial: we consider one-sided generalized
Fourier expansions, the Jacobson radical, invariant subspaces,
spatially implemented automorphisms, reflexivity and Beurling type theorems. Much
of this analysis goes into the proof of our main classification theorem, that the
algebras are unitarily equivalent if and only if their directed graphs are isomorphic.

\vspace{0.1in}

{\noindent}{\it Note added in proof.} This paper has motivated a
number of recent papers on non-selfadjoint graph algebras and
related topics including
\cite{Duncan,JP,JK2,JuryK,KK2,KK1,KP3,KP2,Solel}.

%%%%%%%%%%%%%%%%%%%%%%%%%%%%%%%
\section{Main Results}\label{S:intro}
%%%%%%%%%%%%%%%%%%%%%%%%%%%%%%%

In this section we outline the main results of the paper.
Every countable directed graph $G$ generates in a recursive way a tree graph and an associated
Hilbert space $\H_G$. On this Hilbert space, which can be viewed as a generalized Fock space,
there are natural `partial creation operators'  $\{ L_e \}$, one for each directed edge, and
projections $\{L_x\}$, one
for each vertex \cite{Muhly_fd}. Alternatively, the family $\{ L_e, L_x\}$  arises
through the left regular representation $\lambda_G$ of what we term the {\it free semigroupoid
of $G$} determined by the
directed paths and vertices in $G$. We develop a structure theory for
the $\wot$-closed algebras $\flgee$ generated by families $\{L_e, L_x\}$, which we call {\it free
semigroupoid algebras}.

There is also a natural right regular representation $\rho_G$ determined by $G$. This yields
partial isometries and projections $\{ R_e ,R_x\}$ on $\H_G$ which commute with $\{L_e,
L_x\}$. In fact, the $\wot$-closed algebra $\frgee$ generated by $\{R_e, R_x\}$ coincides with
the commutant $\flgee^\prime = \frgee$. Furthermore, $\frgee^\prime = \flgee$ and hence
$\flgee$ is its own second commutant. We establish this by first observing
that $\flgee$ is unitarily equivalent to $\fR_{G^t}$, where
$G^t$ is the {\it transpose graph} with directed graph obtained from $G$ simply by reversing
the directions of all edges.  In general, the
algebras $\frgee$ and $\fR_{G^t}$ can be
different and so this is a departure from the free semigroup case, where
$\fL_n$ and $\fR_n$ are unitarily equivalent just by symmetry. The methods here provide
Fourier expansions for all elements
of $\flgee$, giving us a key technical device.

The question of semisimplicity is considered in Section~\ref{S:radical} and we
prove  $\flgee$ is semisimple precisely when $G$ is transitive in each component; that is, every
edge lies on a cycle.
Also we show  for a  finite graph $G$ that the Jacobson radical of $\flgee$ is
$\wot$-closed, nilpotent and  equal to
the  ideal generated by the $L_e$ which correspond to edges $e\in\edges(G)$
which do not lie on a cycle in $G$.
Further, we show there is a block matrix decomposition of $\B(\H_G)$ such that  the
radical is the off-diagonal component of $\flgee$ in this decomposition.

In Section~6 we discuss several examples illustrating the diversity of the
algebras $\flgee$. In particular, a matrix function representation is obtained in
the case of the graphs which are a single directed cycle. Also we indicate how
amalgamations of graphs correspond to the ($\wot$-closed) free product with
amalgamation of the associated free semigroupoid algebras.

We begin the analysis of  the invariant subspace structure in
Section~\ref{S:invariant}. We first prove $\flgee$ is reflexive by
considering some obvious invariant subspaces associated with
$\frgee$. Next the eigenvalues and eigenvectors for $\flgee^*$ are
computed. In Section~\ref{S:beurling} we prove an invariant
subspace theorem of Beurling type \cite{Beur,DP1,MS_CJM,Pop_fact}.
Every invariant subspace of $\flgee$ is generated by a wandering
subspace for the algebra, and is the orthogonal direct sum of
cyclic subspaces which are minimal in the sense that the cyclic
vector is supported by some projection $L_x$. Each of these
minimal cyclic invariant subspaces is the range of a partial
isometry in $\frgee$, and the choice of partial isometry is unique
up to a scalar multiple. We then use the Beurling Theorem to
derive an explicit characterization of partial isometries in
$\flgee$, and this yields an {\it inner-outer factorization}
\cite{Douglas_text,DP1,Hoffman_text,Pop_beur} for elements of
$\flgee$.

We resolve the classification problem for these algebras in Section~\ref{S:classification}  by
proving  $G$ to be a complete unitary invariant of $\flgee$. In fact, there are simple
dimension formulae through which one can calculate the graph structure from the
algebra and the ideal $\flgee^0$
generated by the $L_e$. However, this ideal is not necessarily preserved under unitary
equivalence and to correct for this we analyze the unitarily
implemented automorphisms of $\flgee$. In particular, we show that these
automorphisms act transitively on sets of  eigenvectors for $\flgee^*$.

In the final section we introduce the notion of a {\it partly
free} $\wot$-closed operator algebra. Meaning that the algebra
contains a copy of the free semigroup algebra $\fL_2$. In spirit,
this is a non-selfadjoint analogue of the requirement that a
$\ca$-algebra contains a copy of the Cuntz algebra $\O_2$. In the
finite vertex case we prove the algebras $\flgee$ are partly  free
precisely when  $G$ contains a {\it double-cycle}. We also
consider the stronger notion of a  {\it unitally partly free}
$\wot$-closed operator algebra. In this case there is a unital
injection of $\fL_2$ into the algebra, and so it contains a pair
of isometries with mutually orthogonal ranges. For the algebras
$\flgee$ we prove this happens exactly when $G$ satisfies a graph
theoretic condition determined by its double-cycles. We conclude
the paper by discussing hyper-reflexivity for $\flgee$.

%%%%%%%%%%%%%%%%%%%%%%%%%%%%%%%%%%
\section{Free Semigroupoid Algebras}\label{S:freesem}
%%%%%%%%%%%%%%%%%%%%%%%%%%%%%%%%%%

Let $G$ be a finite or countably infinite directed graph with edge
set $E(G)$ and vertex set $V(G) $. Let $\fngee$ be the {\it free
semigroupoid} determined by $G$. By this we mean that  $\fngee$
consists of the set of  vertices $x$ of $G$ and the finite paths
$w$ of edges $e$ in $G$, together with the natural operation of
concatenation of allowable paths. (This is also called the {\it
path space} of $G$.) We view each vertex as a degenerate path.
Given a path $w$ in $\fngee$ we write $w=y wx$ when the initial
and final vertices of $w$ are, respectively, $x$ and $y$. We shall
use the {\it range} and {\it source} maps $r$ and $s$ to indicate
this, so that $s(w)=x$ and $r(w)=y$. At times it will be
convenient to use the {\it transition matrix} $A = (a_{yx})$
associated with $G$, where $a_{yx}$ is the number of directed
edges in $G$ from vertex $x$ to vertex $y$.

We note that the terminology `free' is appropriate as  all paths are admissable elements of
$\fngee$ and there are
no relations amongst words, that is, there are no reducible products in $\fngee$ other than trivial
reductions
involving units. Furthermore, there is clearly a natural  groupoid, an inverse semigroup in fact
(which deserves to be called the free groupoid of $G$)  which contains
$\fngee$ as a semigroupoid, just as in the free group case.

Let $\H_G  = \ell^2(\fngee)$ be the Hilbert space with
orthonormal basis $\{ \xi_w : w\in \fngee\}$ indexed by elements of $\fngee$. For each edge
$e\in E(G)$ and vertex $x \in V(G)$, we may define partial isometries and projections on
$\H_G$ by:
\[
L_e\xi_w = \left\{ \begin{array}{cl}
\xi_{ew} & \mbox{if $r(w)=s(e)$} \\
0 & \mbox{if $r(w)\neq s(e)$}
\end{array}\right.
\]
and
\[
L_x \xi_w = \left\{ \begin{array}{cl}
\xi_{xw}=\xi_w & \mbox{if $r(w)=x$} \\
0 & \mbox{if $r(w)\neq x$}
\end{array}\right.
\]
We shall use the convention $\xi_{ew} = 0$ if $r(w)\neq s(e)$. The
family $\{L_e,L_x\}$ also arises \cite{Muhly_fd} through the left
regular representation $\lambda_G: \fngee \rightarrow \B(\H_G)$,
with $\lambda_G(e) = L_e$, and $\lambda_G(x) = L_x$.

\begin{rem}
{}From the  $\ca$-algebra perspective, the partial isometries $\{L_e\}$ are of
Cuntz-Krieger-Toeplitz type in the sense that the $\ca$-algebra generated by $\{ L_e \}_{e\in
E(G)}$ is generally the extension by the compact operators of a Cuntz-Krieger $\ca$-algebra
\cite{BHRS,Kumjian1,Kumjian2}.

As we see in later examples,
there is a useful  interpretation of the actions of the operators
$L_e$ (and the companion operators $R_e$) in terms of
the  natural tree graph whose nodes are labelled by the paths $w$.
This tree is
generated recursively from the `tree top graph' consisting of the edges of $G$
directed down from the vertices of $G$.
The nodes $w$  correspond to the basis vectors $\xi_w$ and the operators
$L_e$ and $R_e$ correspond to visualizable partial bijections of the tree structure.
The tree perspective is also useful in the analysis of
eigenvectors for $\flgee^*$ in Section~\ref{S:invariant}.
\end{rem}

\begin{defn}
The {\it free semigroupoid algebra} determined by $G$ is the weak operator
topology closed algebra generated by  $\{L_e,L_x\}$,
\begin{eqnarray*}
\flgee &=& \wot - \Alg\,\, \{ L_e,L_x :e\in E(G), x\in V(G) \} \\
&=& \wot - \Alg \,\,\{ \lambda_G(w) : w\in\fngee \}.
\end{eqnarray*}
\end{defn}

Finally we define $\flgee$, up to unitary equivalence, in terms of generators and relations
together with a spatial condition. This leads naturally into the topic of representation theory of
$\flgee$ which we will not pursue in this paper. Nevertheless, we introduce the basic idea of a
free partial isometry representation $\pi$ of a directed graph $G$, which in turn gives rise to an
operator algebra $\fL_\pi$ which one might refer to as a free semigroupoid algebra. In the purely
atomic case, which is also characterized by a simple spatial condition, the algebras $\fL_\pi$ are
all naturally isomorphic to $\flgee$.

\begin{defn}
Let $G$ be a countable directed graph and let $\H$ be a separable Hilbert space. A
{\it free partial isometry representation} $\pi$ of $G$ is a pair of maps
\[
\pi : V(G) \longrightarrow \proj (\H), \,\,\,\,\, \pi : \edges(G) \longrightarrow \pisom(\H),
\]
denoted $x \rightarrow P_x$ and $e \rightarrow S_e$ respectively, such that
\begin{itemize}
\item[$(i)$] given $s(e)=x$, the initial projection $S_e^* S_e =
P_x \neq 0$, \item[$(ii)$] the projections $P_x$ have pairwise
orthogonal ranges and sum to the identity operator $I$,
\item[$(iii)$] the projections $\{ S_e S_e^*: e\in\edges(G)\}$ are
pairwise orthogonal and for each $x\in V(G)$,
\[
E_x = P_x - \sum_{r(e)=x} S_e S_e^* \geq 0.
\]
\end{itemize}
\end{defn}

Such a partial isometry representation of $G$ is indeed free in the sense that for each word $w$
in $\fngee$ the corresponding operator $\pi(w)$ (defined as the natural product) is a non-zero
partial isometry and so we obtain a faithful representation of $\fngee$. We write $\fL_\pi$ for
the weak operator topology closed algebra generated by $\{ \pi(e), \pi(x) : e\in\edges(G), \,\, x\in
V(G)\}$. Let us further define a {\it purely atomic} free partial isometry representation to be one
for which each of the vacuum projections $E_x$ is non-zero, and the set of final projections
$\pi (w) E_x \pi(w)^*$, for $x$ in $V(G)$ and $w\in \fngee$, sum to the identity. Then we have
the following proposition.

\begin{prop}\label{repntheory}
If $\pi$ is a purely atomic free partial isometry representation of the directed graph $G$, then
$\fL_\pi$ is naturally isometrically isomorphic to $\fL_G$ by an isomorphism which is $\wot -
\wot$ continuous. If the vacuum projections $E_x$, $x\in V(G)$, are each of rank one, then
$\fL_\pi$ and $\flgee$ are unitarily equivalent.
\end{prop}

\Prf
If the vacuum projections each have rank one, then it follows that for a
choice of unit vacuum
vectors $\eta_x$ with $E_x \eta_x = \eta_x$, the vectors $\eta_w =\pi(w) \eta_x$, for $wx=w$,
make up an orthonormal basis. The correspondence $\eta_w \rightarrow \xi_w$ gives the desired
unitary equivalence of $\fL_\pi$ and $\fL_G$.
If the vacuum projections are of infinite rank, then $\pi$ is unitarily equivalent to the free partial
isometry representation $\lambda_G^{(\infty)}$. Thus we have the natural
isomorphisms
$\fL_\pi \simeq \fL_\pi^{(\infty)}\simeq \flgee^{(\infty)} \simeq \flgee. $
\bx

%\begin{rem}
% For the operator algebras $\fL_\pi$ given in
%Proposition~\ref{repntheory} it is natural to ask for
%complete unitary equivalence invariants which involve the defect
%multiplicities $\dim \big( E_x
%\big)$. However, to resolve this, and thereby generalize the
%Beurling-Lax-Halmos classification
%of pure isometries, we must first understand unitary automorphisms of the
%algebras $\flgee$
%themselves. This we do in Section~\ref{S:classification}.
%\end{rem}

%\begin{lem}\label{projnlemma}
%For $w\in\fngee$ and all $i$ we have
%Similarly, $L_w P_i=L_w$ if $w=wx_i$, and $L_wP_i=0$ otherwise. Further,
%for all $w\in\fngee$, either $L_w^*L_w =0$ or $L_w^*L_w =P_i$ where
%$w=wx_i$. There are analogous properties for the $R_w$ and $Q_i$.
%\end{lem}

%\Prf
%The action of $P_iL_w$ is described on basis vectors by $P_i L_w \xi_v =
%\xi_{x_iwv}$, hence the above identities are apparent. The behaviour of
%$L_wP_i$ is similar. If $w = e_{i_k}\cdots e_{i_1}$, then
%\[
%L_w^*L_w = L_{e_{i_1}}^* \cdots L_{e_{i_{k-1}}}^*
%(L_{e_{i_k}}^*L_{e_{i_k}})L_{e_{i_{k-1}}} \cdots L_{e_{i_1}}.
%\]
%Thus, since each $L_e^*L_e = P_j$ for some $j$, the previous discussion
%shows that $L_w^*L_w= P_i$ where $w=wx_i$ when this operator is non-zero.
%\bx

%%%%%%%%%%%%%%%%%%%%%%%%%%%%%%%%%%%%%%%%%%%%%%
%%%%%%%%%%%%
\section{Commutant and Basic Properties}\label{S:commutant}
%%%%%%%%%%%%%%%%%%%%%%%%%%%%%%%%%%%%%%%%%%%%%%
%%%%%%%%%%%%

Let $G$ be a directed graph, perhaps with a countable number of
edges or vertices. Given $w = e_{i_k} \cdots e_{i_1} \in \fngee$,
let $L_w$ be the partial isometry $L_w = L_{ e_{i_k}} \cdots L_{
e_{i_1}}$. We shall also use the notation $P_x \equiv L_x$ for
vertex projections, or $P_i \equiv L_{x_i}$ if $V(G) = \{ x_1,
x_2, \ldots \}$ has been enumerated. Observe that $L_w = P_{r(w)}
L_w P_{s(w)}$. Define another family $\{ R_e, R_x \}$ of partial
isometries and projections on $\H_G$ by the equations $R_e \xi_w =
\xi_{we}$ and $R_x \xi_w = \xi_{wx}$. These operators also come
from the right regular representation
$\rho_G:\fngee\longrightarrow \B(\H_G)$, which yields partial
isometries $\rho_G(w) \equiv R_{w^\prime}$ for $w\in\fngee$ acting
on $\H_G$ by the equations $R_{w^\prime}\xi_v = \xi_{vw}$, where
$w^\prime$ is the word $w$ in reverse order. Observe that
$R_{v^\prime} L_w = L_w R_{v^\prime}$ for all $v,w\in\bbF^+(G)$.
In fact, in this section we show the algebra
\begin{eqnarray*}
\frgee &=& \wot - \Alg \,\, \{R_e,R_x :e\in E(G), x\in V(G) \} \\
&=& \wot - \Alg \,\, \{ \rho_G(w) : w\in\fngee \}
\end{eqnarray*}
coincides with the commutant $\flgee^\prime = \frgee$. We shall
use notation $Q_x \equiv R_x$ for the vertex projections in
$\frgee$. Thus $Q_x$ is the projection onto the subspace spanned
by all $\xi_w$ with $s(w) =x$. Let $G^t$ denote the {\it transpose
graph} of $G$. This is the directed graph obtained from $G$ simply
by reversing directions of all directed edges. If $v = e_{i_1}
\cdots e_{i_k}$ is a product of edges in $E(G)$, then we let $v^t$
be the product $v^t = e_{i_k}^t \cdots e_{i_1}^t$, where $e^t$ is
the directed edge in $G^t$ which is $e_i$ with direction reversed.
Further define $x^t = x$ for $x\in V(G) = V(G^t)$.

\begin{lem}\label{commutantlemma}
The algebras $\fL_G$ and $\fR_{G^t}$ are unitarily equivalent via the map
$W: \H_{G^t} \rightarrow \H_G$ defined by $W\xi_{v^t} =
\xi_v$.
\end{lem}

\Prf
The map $W$ is easily seen to be a unitary operator.
Given $v\in \fngee$ and $e\in E(G)$,
\[
(W^* L_e W) \xi_{v^t} = W^* L_e \xi_v = W^* \xi_{ev}
= R_{e^t} \xi_{v^t}.
\]
Hence $W^*L_eW = R_{e^t}$ for $e\in\edges(G)$, and similarly $W^* L_x W = R_{x^t}$ for
$x\in V(G)$.  Thus we have $W^* \flgee W = \fR_{G^t}$.
\bx

As in the free semigroup case \cite{DP1}, we can consider the Cesaro
operators
associated with the partition $I= E_0+E_1+\ldots$ where $E_k$ is the
projection onto $\spn\{\xi_w : |w|=k, w\in \fngee\}$, the subspace spanned by basis vectors from
paths of length $k$. These operators are given by
\[
\Sigma_k(A) = \sum_{|j|<k} \Big( 1-\frac{|j|}{k}\Big) \Phi_j(A),
\]
where the operators $\Phi_j(A) = \sum_{k\geq \max\{ 0,-j\}} E_k A E_{k+j}$
are the diagonals of $A$ with respect to this block
decomposition, and $\Sigma_k(A)$ converges to $A$ in the strong operator
topology for all $A \in\B(\H_G)$.

We mention that as a consequence of  Lemma~\ref{commutantlemma}
one sees that the following result is a generalization of
Proposition~5.4 from \cite{MS_CJM} where the commutant was shown
to be determined by the algebra associated with the transpose
graph.

\begin{thm}\label{commutant}
The commutant of $\fR_G$ coincides with $\flgee$.
\end{thm}

\Prf We have observed above that $\flgee$ is contained in
$\frgee^\prime$. To see the converse, fix $A \in \frgee^\prime$.
We show that $A_x \equiv AP_x$ belongs to $\flgee$ for all $x\in
V(G)$. Let $A\xi_x =Q_x A_x\xi_x= \sum_{s(w)=x} a_w \xi_w.$ Define
operators in $\flgee$ by
\[
p_k(A_x) = \sum_{ |w|<k; \, s(w)=x} \Big( 1 - \frac{|w|}{k}\Big)
a_w L_w.
\]
We claim that $A_x = \sotlim_{k\rightarrow\infty} p_k(A_x)$. This
will be proved by showing that $p_k(A_x) =\Sigma_k(A_x)$. First
observe that $\Phi_j(A_x)\in\frgee^\prime$ for all $j$. Indeed,
this operator commutes with the $R_e$ since $A_x$ belongs to
$\frgee^\prime$ and $E_{k+1} R_e = R_e E_k$ for all $k$, while
$\Phi_j(A_x)$ commutes with each $Q_y$ since $Q_y E_k = E_k Q_y$
is the projection onto $\spn\{ \xi_w: |w| =k, \,\, s(w)=y \}$ for
all $k$. It follows that $\Sigma_k(A_x)$ belongs to
$\frgee^\prime$ for $k\geq 1$.

Now it is enough to show that
$\Sigma_k(A_x) \xi_x = p_k(A_x) \xi_x$. If this is the case, then for
$w=xw$ in $\fngee$
\[
\Sigma_k(A_x) \xi_w = R_w \Sigma_k(A_x) \xi_x = R_w p_k(A_x) \xi_x =
p_k(A_x) \xi_w,
\]
Whereas for $w=yw$ with $y\neq x$ we have
\[
\Sigma_k(A_x)\xi_w = \Sigma_k(A_x)P_x P_y\xi_w=0= p_k(A_x)P_y \xi_w =
p_k(A_x)\xi_w,
\]
since $P_x$ commutes with each $E_l$.
However, observe that $\Phi_0 (A_x) \xi_x = E_0(A_x)E_0\xi_x = a_x
\xi_x$, and $\Phi_j(A_x)\xi_x = 0$ for $j>0$. Further, for $j< 0$ we have
\[
\Phi_j(A_x) \xi_x = (E_{-j} A_x )\xi_x = E_{-j} \sum_{s(w)=x}
a_w\xi_w = \sum_{s(w)=x; \,|w|=-j} a_w \xi_w.
\]
Hence it follows that
\[
\Sigma_k(A_x) \xi_x =   \sum_{|w| <k;\, s(w)=x} \Big( 1
-\frac{|w|}{k}\Big) a_w \xi_w = p_k(A_x) \xi_x.
\]

We have established that each $A_x = AP_x$ belongs to $\flgee$.
This completes the proof since $A=\sum_{x\in V(G)} AP_x$, the sum
converging in the strong operator topology when $V(G)$ is
infinite. \bx

\begin{rem}\label{Fourier}
{} From the proof of this theorem elements of $\flgee$ can be seen
to have Fourier expansions. In particular, if $A$ belongs to
$\flgee$ with $A\xi_x= AP_x\xi_x = Q_x (AP_x)\xi_x = \sum_{s(w)=x}
a_w \xi_w$ for $x\in V(G)$, then $A\xi_v = R_v A\xi_x =
\sum_{s(w)=x} a_w \xi_{wv}$ for $v=xv \in \fngee$, and it follows
that the Cesaro partial sums associated with the series
$
A \sim  \sum_{w\in \bbF^+(G)} a_w L_w
$
converge in the strong operator topology to $A$.
\end{rem}

\begin{cor}\label{leftcommutant}
The commutant of $\flgee$ coincides with $\frgee$.
\end{cor}

\Prf
If $W$ is the unitary from Lemma~\ref{commutantlemma}, we have
\[
\fR_{G^t}^\prime = (W^*\flgee W)^\prime = W^* \flgee^\prime W.
\]
On the other hand, it also follows from the lemma and the
definition of $W$ that $\fR_G = W\fL_{G^t}W^*$. But the theorem
tells us $\fL_{G^t} = \fR_{G^t}^\prime$ and hence  $ \fR_G = W
\fR_{G^t}^\prime W^* = \flgee^\prime. $ \bx

The following are simple consequences of the previous two results.

\begin{cor}
Let $G$ be a countable directed graph.
\begin{itemize}
\item[$(i)$] $\flgee$ is its own second commutant, $\flgee
=\flgee^{\prime\prime}$.
\item[$(ii)$] $\flgee$ is inverse closed.
\end{itemize}
\end{cor}

We can also describe the self-adjoint part of these algebras.

\begin{cor}\label{normals}
The set of normal elements in $\flgee$ is precisely  the  span of
$\{ P_x : x\in V(G) \}$. Similarly, the normal elements in $\frgee$ belong to the  span
of $\{ Q_x : x\in V(G) \}$.
\end{cor}

\Prf The operators belonging to the  span of the $P_x$ are normal
elements since the $P_x$ are projections with pairwise orthogonal
ranges. Let $A$ be a normal element of $\flgee$ and put $\alpha_x
= (A \xi_x, \xi_x)$ for $x\in V(G)$. Clearly $\xi_x$ is an
eigenvector for $\flgee^*$, since it is an eigenvector for each of
the generators, hence $A^* \xi_x = \ol{\alpha}_x \xi_x$ and by
normality $A \xi_x = \alpha_x \xi_x$. However, recall that $ \ran
(P_x) = \spn\{\xi_w : r(w)=x  \} = \spn\{ R_w \xi_x : r(w)=x \}. $
Thus, as $A$ commutes with $\fR_G$ we have
\[
A \xi_w = R_w A \xi_x = \alpha_x \xi_w \qfor w=xw\in\fngee.
\]
In other words, $A P_x = \alpha_x P_x$ and it follows that $A =
\sum_{x\in V(G)} \alpha_x P_x$.
\bx

\section{The Radical}\label{S:radical}
%%%%%%%%%%%%%%%%%%%%%%%%%%%%%%%%%%%%%%%%%%%%%%
%%%%%%%%%

In this section we determine when $\flgee$ is semisimple and in
the case of finite graphs we show how to compute the Jacobson
radical of $\flgee$ strictly in terms of the graph structure.  In
particular, there is a block matrix decomposition of $\flgee$ in
which the radical is the off-diagonal part.

A directed  graph $G$ is {\it transitive} if there are paths in
both directions between every pair of vertices in $G$.  A {\it
(connected) component} of $G$ is given by a maximal collection of
vertices and edges which are joined in the undirected graph
determined by  $G$.
%\Prf
%First note that if $e = xey$ is an edge of a transitive component of $G$, then
%there is a path $u =
%yux$, and hence $eu = x(eu)x$ is a cycle. On the other hand, if $x,y$ are
%vertices in the same
%connected component, then there is an undirected path between $x$ and $y$. If $G$
%is a cycle
%graph, it follows that there is also a directed path between $x$ and $y$ and vice
%versa (This
%point is clarified in the proof of Lemma~\ref{cyclelemma}).
%\bx
By a {\it cycle}, we mean a path in $G$ with the same initial and
final vertices. It is not hard to see that $G$ is transitive in
each component precisely when every directed edge in $G$ lies on a
cycle. Let $\badedges(G) \subseteq E(G)$ be the collection of
edges $e\in\edges(G)$ which {\it do not} lie on a  cycle. The set
$\badedges(G)$ is empty precisely when $G$ is transitive in each
component. The Jacobson radical is determined by these edges.

\begin{thm}\label{semisimple}
$\flgee$ is semisimple if and only if $G$ is transitive in each
component. When  $G$ has finitely many vertices, $|V(G)|= M <
\infty$, the radical  is nilpotent of degree at most $M $ and is
equal to the $\wot$-closed two-sided ideal generated by $\{ L_e :
e\in\badedges(G)\}$.
\end{thm}

%Let us begin by observing a feature of paths which avoid $\badedges(G)$.

%\begin{lem}\label{cyclelemma}
%For every path $v\in\fngee$ which consists of
%edges not belonging to $\badedges(G)$, there is a path $u\in\fngee$ such
%that $uv\in\fngee$ forms a  cycle.
%\end{lem}

%\Prf
%For the sake of brevity, we shall prove the lemma for $v=fe$ with each
%edge $e,f$ part of a cycle in the directed graph for $G$. The general case
%follows in a similar manner. The path $v$ moves from the initial vertex $x$ of
%$e$ along the directed edge $e$ to its final vertex $y$, which is the
%initial vertex of $f$, then follows the edge $f$ to its final vertex
%$z$. As $f$ lies on a cycle we can find a path $u_1=y u_1 z$, and
%since
%$e$ lies on a cycle we can find a path $u_2 = x u_2 y$. Thus we have
%the
%allowable path
%\[
%(x u_2 y)(y u_1 z)(z f y)(y e x) = x (u_2u_1fe) x = x (u_2u_1 v)x,
%\]
%which is a cycle containing $v$, as required.
%\bx

We begin by proving one of the implications in the theorem.

\begin{lem}\label{sslemma1}
If $G$ is transitive in each component, then $\flgee$ is semisimple. In particular,  for every
non-zero $A$ in
$\flgee$, there is a path $w\in\fngee$ such that $L_w A$ is
not quasinilpotent.
\end{lem}

\Prf By Theorem~\ref{commutant}, the Fourier expansion of
$A\in\flgee$ is determined by the vectors $\{ A\xi_x : x\in
V(G)\}$. If $A$ is non-zero, then some $A\xi_x = \sum_{s(w)=x}
a_{w} \xi_w$ is non-zero and there is   a path $v$ of minimal
length such that $a_{v}\neq 0$. As $G$ is transitive in each
component, there is a path $u\in\fngee$ such that $uv$ is a cycle.
Hence the paths $(uv)^k$ for $k\geq 1$ are cycles in $\fngee$.
Thus for $k\geq 1$, the expansion for $L_uA$ gives us
\[
(L_uA)^k \xi_x = (L_u A)^{k-1} \sum_{s(w)=x} a_{w} \xi_{uw} =
a^k_{v} \xi_{(uv)^k} + \sum_{w\neq (uv)^k} b_w \xi_w.
\]
Thus for $k\geq 1$ we have
\[
||(L_uA)^k||^{1 / k} \geq \big|( (L_uA)^k\xi_x,
\xi_{(uv)^k})\big|^{1 / k} = |a_{v}| > 0,
\]
so the operator $L_uA$ has positive spectral radius and is not
quasinilpotent. But recall the radical $\rad \flgee$ is equal to
the largest quasinilpotent ideal in $\flgee$. Thus $\rad \flgee =
\{ 0 \}$ when $G$ is transitive in each component. \bx

Towards the converse implication we  obtain the following
description of edges not lying on cycles.

\begin{lem}\label{sslemma2}
The following assertions are equivalent for
$e\in\edges(G)$.
\begin{itemize}
\item[$(i)$] $L_e\in\rad \flgee$.
\item[$(ii)$]  $e\in\badedges(G)$.
\item[$(iii)$] $(AL_e)^2=0$ for all $A\in\flgee$.
\end{itemize}
\end{lem}

\Prf By considering Fourier expansions, the last two conditions
are easily seen to be equivalent to the requirement $L_w^2 =
L_{w^2} =0$ whenever $w\in \fngee$ is a path which includes $e$ as
an edge. The characterization of $\rad \flgee$ used in the
previous proof shows that $(iii)\Rightarrow (i)$. Finally, if
$(ii)$ fails, then there is a path $u\in \fngee$ such that
$(ue)^k$ is a cycle for $k\geq 1$, and we may argue as in the
previous proof to show that $(i)$ fails. \bx

\vspace{0.08in}

{\noindent}{\it Proof of Theorem~\ref{semisimple}.}
The algebra $\flgee$ is semisimple when $G$ is transitive in each component by
Lemma~\ref{sslemma1}. On the other hand, if there is a component in $G$ which is not
transitive, then the set $\badedges (G)$ is nonempty and Lemma~\ref{sslemma2} gives an edge
$e\in\badedges(G)$ with $L_e\in \rad \flgee$. Thus $\flgee$ has non-zero radical in this case.
It remains to show the radical is nilpotent of degree at most $M$ and equal to the $\wot$-closed
ideal generated by $\{ L_e : e \in B(G)\}$ when the cardinality of the vertex set $V(G)$ is  $M <
\infty$.

Let $\J$ be the $\wot$-closed two-sided ideal in $\flgee$ generated by
$\{L_e : e\in \badedges(G)\}$. We first observe that the radical contains
this ideal. For this we may use the following block matrix decomposition of $\flgee$.
We say that  a subset $H$ of edges and vertices of $G$ is  {\it maximally transitive}  if: there
are directed paths in both directions between every pair of vertices in $H$; the initial and final
vertices of every edge in $H$ belong to $H$;  each edge between every pair of vertices in $H$
also belongs to $H$; and $H$ is maximal with respect to these properties.
Let $\{ G_i \}_{i\in\I}$ be the maximally transitive components of $G$, and let
$\{S_i\}_{i\in\I}$ be the
projections $S_i = \sum_{x\in V(G_i)} P_x$. Then $I = \big( \sum_{i\in\I} \oplus S_i \big)
\oplus \big( \sum_{x\notin \cup_{i\in\I} V(G_i)} \oplus P_x \big)$ and we may
consider the block matrix form of $\flgee$ with respect to this spatial decomposition.
By considering Fourier expansions for elements of $\flgee$, it is not hard to see that the ideal
$\J$ is given by the off-diagonal entries of $\flgee$ in this
decomposition. It follows that $\J^M = \{0\}$, and for all $X\in\flgee$ and $A\in\J$ we have
$(XA)^M =
0$. Hence $\J$ is contained in $\rad \flgee$ and is nilpotent of degree at most $M$.

For the converse inclusion, suppose $A$ belongs to $\rad \flgee$ with expansion
 scalars $\{a_w\}_{w\in\fngee}$. We claim that a
coefficient $a_w$ is non-zero only if the path $w$ includes an edge $e\in\badedges(G)$. This
will complete the proof, since the Cesaro sums for $A$ would then belong to
$\J$, and they converge in the strong operator topology to $A$.
Suppose by way of contradiction that there is a path $v$ with $a_v\neq 0$ which
includes no edges from $\badedges(G)$, and assume $v$ is
a path of minimal length with this property.
Since every path in $G$ which includes no edge from
$\badedges(G)$ must be part of some cycle, we may choose a path $u\in\fngee$ such that $uv =
x( uv) x$ is a cycle in $G$.  Hence $(uv)^k$ belongs to $\fngee$ for
$k\geq 1$. But it is clear that $(uv)^M$ is a path of minimal length amongst the paths
in the expansion of $(L_u A)^M$ with non-zero coefficients. Further, the coefficient of
$L_{(uv)^M}$ in this expansion is $(a_v)^M$.
Thus we have
$
||(L_uA)^{Mk}||^{1 / k} \geq | a_v|^M > 0$ for $k\geq 1.$
Hence $(L_uA)^M = ((L_uA)^{M-1} L_u)A$ has positive spectral radius and is not
quasinilpotent. This is a contradiction since $A$ belongs to $\rad \flgee$, and the result follows.
\bx

\begin{rem}
{} From the block matrix decomposition of $\flgee$ used in the proof of
Theorem~\ref{semisimple} for the finite vertex case, we see the ideal $\rad \flgee$ is given
by the $\wot$-closed ideal determined by the off-diagonal entries in this decomposition, which
in turn is the $\wot$-closed two-sided
ideal generated by $\{ L_e : e\in \badedges(G)\}$. This point is discussed further in the next
section.
We also mention that the general ideal structure of $\flgee$ has been characterized in
\cite{JuryK}.
\end{rem}

%%%%%%%%%%%%%%%%%%%%%%%%%%%%%%%%%%%%%%%%%
\section{Examples and Amalgamated Free Products}\label{S:freeprods}
%%%%%%%%%%%%%%%%%%%%%%%%%%%%%%%%%%%%%%%%%

We now examine several simple examples of graphs $G$ and their operator algebras $\flgee$.
Some of these algebras admit natural representations as algebras of matrices or as subalgebras of
(possibly infinite) matrix algebras over $H^\infty$. In this case some of their algebraic features
such as the Jacobson radical and the structure of ideals become more apparent. The examples
will also be used to illustrate the discussion at the end of this section where we consider the
structure of $\flgee$ as an amalgamated free product.

\begin{eg}\label{singlevertexeg}
We first consider the single vertex cases.
The algebra generated by the graph with a single vertex $x$ and single loop edge $e = xex$ is
unitarily
equivalent to the classical analytic Toeplitz algebra $H^\infty$
\cite{Douglas_text,Hoffman_text}. Indeed, $\H_G$ in this case may be naturally
identified with the Hardy space $H^2$, and under this identification $L_e$ is easily seen to be
unitarily equivalent to the unilateral shift $U_+$.

The noncommutative analytic Toeplitz algebras $\fL_n$, $n\geq 2$
\cite{AP,DKP,DP1,DP2,Kribs_factor,Pop_fact,Pop_beur}, the
fundamental {\it free semigroup algebras}, arise from the graphs
with a single vertex and $n$ distinct loop edges. For instance, in
the case $n =2$ with loop edges $e=xex \neq f=xfx$, the space
$\H_G$ is identified with unrestricted 2-variable Fock space
$\H_2$. Under this identification, the operators $L_e, L_f$ are
unitarily equivalent to the natural creation operators on $\H_2$
which are the canonical  Cuntz-Toeplitz isometries. Further, $P_x
= I$, and thus $\flgee \simeq \fL_2$.
\end{eg}

\begin{eg}\label{finitedimeg}
If $G$ is a finite directed graph with no directed cycles, then the
Fock space $\H_G$ is finite dimensional and so too is $\flgee$. For example,
consider the graph with three vertices and two edges, labelled $x_1$, $x_2$, $x_3$, $e$, $f$
where $e = x_2ex_1$, $f= x_3fx_1$. Then the Fock space is spanned by the vectors $\{
\xi_{x_1}, \xi_{x_2}, \xi_{x_3}, \xi_e, \xi_f \}$ and a little reflection reveals that with this basis
the general operator
$
\alpha L_{x_1} + \beta L_{x_2} + \gamma L_{x_3} + \lambda L_e + \mu L_f
$
in $\flgee$ is represented by the matrix
\[
\begin{bmatrix}
\alpha & & & & \\
 & \beta & & & \\
 & & \gamma & & \\
\lambda & & & \beta & \\
\mu & & & & \gamma
\end{bmatrix}
\]
(unmarked entries are zero). Alternatively one can reorder the basis to view $\flgee$ as a
subalgebra of $\M_5(\bbC)$ consisting of matrices of the form
\[
\begin{bmatrix}
\alpha & 0 \\ \lambda & \beta
\end{bmatrix}
\oplus
\begin{bmatrix}
\alpha & 0 \\ \mu & \gamma
\end{bmatrix}
\oplus
\begin{bmatrix}
\beta
\end{bmatrix}
\oplus
\begin{bmatrix}
\gamma
\end{bmatrix}.
\]
Algebraically, $\flgee$ is isometrically isomorphic to the so-called digraph algebra in
$\M_3(\bbC)$ consisting of the matrices
\[
\begin{bmatrix}
\alpha & 0 & 0 \\ \lambda & \beta & 0 \\ \mu & 0 & \gamma
\end{bmatrix}.
\]

Recall that a digraph algebra $\A(H)$ is a unital subalgebra of $\M_n(\bbC)$ which is spanned
by some of the standard matrix units of $\M_n(\bbC)$. The graph $H$ is transitive and reflexive
and is such that the edges of $H$ naturally label the relevant matrix units. If we take $H$ to be
the augmentation of our graph $G$ by loops at the vertices, then we can view $\flgee$ as a
faithful representation of $\A(H)$. Notably this representation is not a star extendible
isomorphism since $\A(H)$ and $\flgee$ generate different $\ca$-algebras. In general, for a
finite cycle-less graph $G$, $\flgee$ is isometrically isomorphic to $\A(\td{G})$ where
$\td{G}$ is the transitive completion of $G$ with vertex loops added.

The commutant of $\flgee$ is best understood through Theorem~\ref{commutant}. However,
one can confirm directly that the commutant of $\flgee$ for this example consists of the matrices
\[
\begin{bmatrix}
a & & & & \\
 & b & & & \\
 & & c & & \\
 & \lambda  & & a & \\
 & & \mu & & a
\end{bmatrix}
\]
with $a,b,c,\lambda,\mu$ in $\bbC$. These matrices in fact correspond to the operators
$
a R_{x_1} + b  R_{x_2} + c  R_{x_3} + \lambda R_e + \mu R_f
$
which are the typical elements of $\frgee$.

In the terminology of the second author \cite{Power_approx} the operator algebras $\flgee$ are
finitely acting, since they act on finite dimensional Hilbert spaces. (This is a stronger notion than
finite dimensionality for an operator algebra.) Also the discussion above shows that the operator
algebras $\flgee$ are the so-called inflation algebras of digraph algebras given in
\cite{Power_approx}.
\end{eg}

\begin{eg}\label{easymtxfnalg}
For  a simple matrix function algebra, we may consider the graph $G$ with
vertices $x, y$ and edges $e=xex$, $f=yex$. Then $\flgee$ is generated by $\{ L_e, L_f , P_x ,
P_y \}$. If we make the natural identifications $\H_G = P_x \H_G \oplus P_y \H_G \simeq H^2
\oplus H^2$ (respecting word length), then
\[
L_e \simeq \left[
\begin{matrix}
U_+ & 0 \\
0 & 0
\end{matrix} \right]
\quad
L_f \simeq \left[
\begin{matrix}
0 & 0 \\
U_+ & 0
\end{matrix} \right]
\quad
P_x \simeq \left[
\begin{matrix}
I & 0 \\
0 & 0
\end{matrix} \right]
\quad
P_y \simeq \left[
\begin{matrix}
0 & 0 \\
0 & I
\end{matrix} \right].
\]
Thus, $\flgee$ is unitarily equivalent to
\[
\flgee \simeq \left[
\begin{matrix}
H^\infty & 0 \\
H^\infty_0 & \bbC I
\end{matrix}\right]
\]
where $H^\infty_0$ is the subalgebra of $H^\infty$ functions $h$
with $h(0) = 0$. With this representation of $\flgee$, it is clear
that $\rad \flgee$ is nilpotent of degree 2 and is given by the
$\wot$-closed ideal generated by $L_f$. The commutant structure is
less evident in this representation, nonetheless, we know it is
generated by $\{ R_e, R_f, Q_x, Q_y\}$ where $Q_y = \xi_y \otimes
\xi_y^*$, $Q_x = I - Q_y$, $R_f = \xi_f \otimes \xi_y^*$ and $R_e
= R_e P_x = P_x R_e$.

By simply adding a directed edge $g = xgy$ to the previous graph
we obtain a very different algebra $\fL_{G^\prime}$. Indeed,
$\fL_{G^\prime}$ is unitarily equivalent to its commutant
$\fL_{G^\prime}^\prime = \fR_{G^\prime} \simeq \fL_{(G^\prime)^t}$
since $(G^\prime)^t \simeq G^\prime$. Furthermore,
$\fL_{G^\prime}$ is  unitally partly free  in the sense of
Section~\ref{S:partfree} because it contains isometries with
mutually orthogonal ranges, for instance, $U = L_e^2 + L_f L_g$
and $V = L_e L_g + L_g L_e$ are isometries which satisfy $U^* V
=0$. The algebra $\fL_{G^\prime}$ will be discussed further in the
context of amalgamated free products below.
\end{eg}

\begin{eg}\label{matrixfntheory}
For an example of an infinite matrix function algebra,
let $G$ be the directed graph with transition matrix $A =\begin{sbmatrix}
1 &1 \\ 0& 1 \end{sbmatrix}$, and
let $e=xex$, $f = y fx$ and $g = y gy$ be the directed edges of $G$. Notice that this graph is
obtained from the graph $G$ of the  previous example by adding the edge $g$. Then
$\flgee$ is generated by three partial isometries $L_e$, $L_f$, $L_g$ and their initial projections
$P_1= L_e^*L_e = L_f^*L_f$ and $P_2 = L_g^* L_g$. These operators act on the Fock space
$\H_G$ with basis indexed by the following tree:
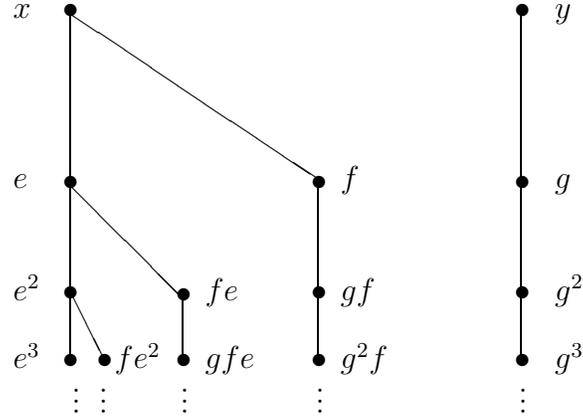
\[
\setlength{\unitlength}{0.15cm}
\begin{picture}(25,38)
%\put(7.6,22){\line(-1,-4){2.6}}
%\put(7.5,22){\line(1,-4){2.6}}
%\put(4.9,21){$x$}
%\put(2.9,11){$e$}
%\put(11,11){$f$}
%\put(6.9,21){$\bullet$}
%\put(4.4,11){$\bullet$}
%\put(9.4,11){$\bullet$}
%\put(1,4){$\bullet$}
%\put(6.7,4){$\bullet$}
%\put(10,12){\line(0,-1){7}}
%\put(9.4,4){$\bullet$}
%\put(11,4){$gf$}
%\put(5,11.5){\line(1,-3){2.2}}
%\put(3.5,4){$fe$}
%\put(5,11.5){\line(-1,-2){3.4}}
%\put(-1.4,4){$e^2$}

\put(31.4,20.8){$\bullet$}
\put(31.4,36){$\bullet$}
\put(35,36){$y$}
\put(35,21){$g$}
\put(32,37){\line(0,-1){15}}
\put(32,21){\line(0,-1){10}}
\put(32,11){\line(0,-1){5}}
\put(31.4,11){$\bullet$}
\put(35,11){$g^2$}
\put(31.4,5){$\bullet$}
\put(35,5){$g^3$}

\put(31.8,1){$\vdots$}
%\put(7.5,-0.5){$\vdots$}

\put(13.4,20.8){$\bullet$}
\put(16,21){$f$}
\put(14,21){\line(0,-1){10}}
\put(14,11){\line(0,-1){5}}
\put(13.4,11){$\bullet$}
\put(16,11){$gf$}
\put(13.4,5){$\bullet$}
\put(16,5){$g^2f$}

\put(13.8,1){$\vdots$}

\put(-8.6,20.8){$\bullet$}
\put(-8.6,36){$\bullet$}
\put(-13,36){$x$}
\put(-13,21){$e$}
\put(-8,37){\line(0,-1){15}}
\put(-8,21){\line(0,-1){10}}
\put(-8,11){\line(0,-1){5}}
\put(-8.6,11){$\bullet$}
\put(-13,11){$e^2$}
\put(-8.6,5){$\bullet$}
\put(-13,5){$e^3$}

\put(-7.8,1){$\vdots$}

\put(2,11){\line(0,-1){5}}
\put(1.4,10.8){$\bullet$}
\put(4,11){$fe$}
\put(1.4,5){$\bullet$}
\put(4,5){$gfe$}

\put(1.8,1){$\vdots$}

\put(-5.6,5){$\bullet$}
\put(-4,5){$fe^2$}
\put(-5.4,1){$\vdots$}

\put(-7.8,36.3){\line(3,-2){22}}

\put(-7.8,21){\line(1,-1){10}}
\put(-7.7,11.3){\line(1,-2){3}}

\end{picture}
\]

Identify the Hardy space $H^2$ of the unit disc with each of the
following `diagonal spaces'; $ \H_1 = \H_e = \spn \big\{ \xi_{x},
\xi_e, \xi_{e^2}, \ldots \big\}, $
\[
\H_2 = L_f \H_e, \,\,\,
\H_3 = L_g L_f \H_e, \,\,\, \ldots \,\,\, ,\H_n = L_g^{n-2} L_f \H_e,\,\,\, \ldots
\]
Also identify the space $\H_g$, similarly defined, with $H^2$. With respect to the
decomposition
$
\H = \Big( \oplus_{n\geq 1} \H_n \Big) \oplus \H_g,
$
the operator $\lambda L_e + \mu L_f + \nu L_g$ has block matrix form
\[
\begin{bmatrix}
        \lambda T_z &0& & &  \\ \mu I & 0 &0& &  \\ &\nu I &0&0& \\ & &\nu I &\ddots &
\ddots  \\
 & & &   \ddots &
       \end{bmatrix} \quad
\oplus \begin{bmatrix}
        \nu T_z
       \end{bmatrix} \quad,
\]
while $\alpha L_e^*L_e + \beta L_g^* L_g = \alpha P_x + \beta P_y$ has the form
\[
\begin{bmatrix}
        \alpha I & &  &  \\  & \beta I & & \\ & &\beta I & \\  & & & \ddots
       \end{bmatrix} \quad
\oplus \begin{bmatrix}
        \beta I
       \end{bmatrix} \quad.
\]
It follows readily that $\flgee$ is naturally unitarily equivalent to the operator algebra of matrix
functions
\[
\begin{bmatrix}
        h_1(z) & &  & &  \\ h_2(z) & \hat{h} (0) &  & & \\ h_3(z) & \hat{h}(1) &\hat{h}(0) & & \\
h_4(z) & \hat{h} (2) & \hat{h} (1) & \hat{h} (0) & \\ \vdots & \vdots & \vdots & & \ddots
       \end{bmatrix} \quad
\oplus \begin{bmatrix}
        h(z)
       \end{bmatrix} \quad
\]
where $h_k \in H^\infty$ with $\sum_{k\geq 1} ||h_k ||^2$ finite, and where $h\in H^\infty$ with
Fourier coefficients $\hat{h}(k)$.

With such an explicit matrix representation, one can examine the ideal structure and other
algebraic aspects in a direct manner. It is clear, for example, that the Jacobson radical of
$\flgee$ is given by the subspace for which $h_1 = h = 0$, and that the quotient by the radical is
isomorphic to $H^\infty \oplus H^\infty$. Less evident is the structure of invariant subspaces for
$\flgee$ which, we shall see, are generated by partial isometries in the commutant of $\flgee$.
\end{eg}

\begin{eg}\label{cyclegraphs}
 The following  algebras play a role in the analysis of
Section~\ref{S:partfree}.
Consider the {\it cycle  graph} $C_n$ which has $n$
vertices $\{x_1, \ldots, x_n \}$  and $n$ edges
$
\big\{ e_1 =x_2 e_1 x_1, \ldots, e_{n-1} =x_n e_{n-1} x_{n-1}, e_n =x_1 e_n x_n \big\}.
$
%\vspace{0.06in}
%\[
%\setlength{\unitlength}{0.15cm}
%\begin{picture}(55,11)
%\put(10,0){\line(0,1){10}}
%\put(20,0){\line(0,1){10}}
%\put(45,0){\line(0,1){10}}
%\put(11.5,9.5){1}
%\put(11.5,-1){2}
%\put(6.5,5){$e_1$}
%\put(21.5,9.8){2}
%\put(21.5,-1){3}
%\put(16.5,5){$e_2$}
%\put(46.5,9.8){N}
%\put(46.5,-1){1}
%\put(40.5,5){$e_N$}
%\put(9.3,-1){$\bullet$}
%\put(9.3,9.5){$\bullet$}
%\put(19.3,-1){$\bullet$}
%\put(19.3,9.8){$\bullet$}
%\put(44.3,-1){$\bullet$}
%\put(44.3,9.8){$\bullet$}
%\put(30,5){$\cdot$}
%\put(33,5){$\cdot$}
%\put(27,5){$\cdot$}
%\end{picture}
%\]
%\vspace{0.13in}
Let $\H \equiv \H_{C_n} = \H_1 \oplus \ldots \oplus \H_n$ be the decomposition of Fock space
corresponding to
the
tree components ($n$ infinite stalks in this case), so that $\H_i = Q_i \H \equiv Q_{x_i} \H$ for
each $i$, and let $\H_{i,k} = L_{e_k}^* L_{e_k}
\H_i = P_k Q_i \H$, for $1 \leq i,k \leq n$. Thus each subspace $\H_i$ breaks up into the direct
sum of $n$ subspaces $\{ H_{i,k}:1\leq k \leq n \}$. With the natural (top down) basis ordering
we have
the identification of each space $\H_{i,k}$ with $H^2$. With this identification note that the
operator
\[
L_{e_i}\,\, : \,\, \H_{i,k} \rightarrow \H_{i,k+1}
\]
(with $\H_{i,n+1} = \H_{i,1}$) is the identity operator from $H^2$ to $H^2$ unless $k = i+1
\,\, ({\rm mod}\, n)$ in which case the operator is the unilateral shift $T_z$. It follows that
$\fL_{C_n}$ is realized as a subalgebra of the $n$-fold direct sum of matrix algebras
$
\M_n (H^\infty) \oplus \ldots \oplus \M_n (H^\infty).
$

For example, when $n = 3$, with respect to the spatial decomposition
\[
\H = \sum_{i=1}^3 \oplus Q_i \H = \left( \sum_{k=1}^3 \oplus \H_{1,k}\right) \oplus
\left( \sum_{k=1}^3 \oplus \H_{2,k}\right) \oplus \left( \sum_{k=1}^3 \oplus \H_{3,k}\right)
\]
the operator $\lambda L_{e_1} + \mu L_{e_2} + \nu L_{e_3}$ has
operator matrix of the form
\[
\begin{bmatrix}
        0 &0&\nu T_z \\ \lambda I & 0 &0 \\0&\mu I & 0
       \end{bmatrix} \quad
\! \oplus \,\,
\begin{bmatrix}
        0 & 0&\nu I \\ \lambda T_z & 0 &0 \\ 0& \mu I & 0
       \end{bmatrix} \quad
\! \oplus \,\,
\begin{bmatrix}
        0 &0 & \nu I \\ \lambda I & 0 &0 \\0& \mu T_z & 0
       \end{bmatrix} \quad
\]
while $\alpha L_{e_1}^* L_{e_1} + \beta L_{e_2}^* L_{e_2} + \gamma L_{e_3}^* L_{e_3}
=\alpha P_1 + \beta P_2 + \gamma P_3$ has the form
\[
\begin{bmatrix}
        \alpha I &0& 0 \\ 0 & \beta I &0 \\0& 0 & \gamma I
       \end{bmatrix} \quad
\! \oplus \,\,
\begin{bmatrix}
        \alpha I & 0& 0 \\ 0 & \beta I &0 \\ 0& 0  & \gamma I
       \end{bmatrix} \quad
\! \oplus \,\,
\begin{bmatrix}
        \alpha I &0 & 0 \\ 0 & \beta I &0 \\0& 0  & \gamma I
       \end{bmatrix}. \quad
\]
It follows that $\fL_{C_3}$ can be identified with the matrix function algebra
\[
\begin{bmatrix}
        h_{11}&z h_{12}& z h_{13} \\ h_{21} & h_{22} & z h_{23} \\  h_{31}& h_{32} &
h_{33}
       \end{bmatrix} \quad
\! \oplus \,\,
\begin{bmatrix}
        h_{11} &  h_{12}& h_{13} \\ z h_{21} & h_{22} &z h_{23} \\ z h_{31}& h_{32} &
h_{33}
       \end{bmatrix} \quad
\! \oplus \,\,
\begin{bmatrix}
        h_{11} & z h_{12}& h_{13} \\  h_{21} & h_{22} & h_{23} \\ z h_{31}& z h_{32} &
h_{33}
       \end{bmatrix} \quad
\]
where $h_{i_j} \in H^\infty$ for $1 \leq i,j \leq 3$.

Once again, one can consider this explicit matrix function algebra representation in the analysis
of the ideals of $\fL_{C_n}$. For example there is a unique maximal ideal, whose intersection
with the centre of $\fL_{C_n}$ identifies with $z H^\infty$, such that the quotient is not
isomorphic to $\M_3 (\bbC)$. This exceptional quotient is isomorphic to the matrix algebra in
$\M_3(\bbC) \oplus \M_3(\bbC) \oplus \M_3(\bbC)$ consisting of the scalar matrices
\[
\begin{bmatrix}
        a &0& 0 \\ d & b &0 \\ f & e & c
       \end{bmatrix} \quad
\oplus \,\,
\begin{bmatrix}
        a & h& i \\ 0 & b  &0 \\ 0& e  & c
       \end{bmatrix} \quad
\oplus \,\,
\begin{bmatrix}
        a &0 & i \\ d & b &j \\0& 0  & c
       \end{bmatrix}. \quad
\]
Also one can verify readily that $\flgee$ is semisimple.

There is an alternative more succinct identification of the {\it cycle algebras} $\fL_{C_n}$
which makes a
connection with semicrossed product algebras. To see this
identify $L_{x_i}\H$ with $H^2$ for each $i$ in the natural way
(respecting word length). Then $\H = L_{x_1}\H \oplus \ldots \oplus L_{x_n}\H \simeq
\bbC^n \otimes H^2$ and the operator $\alpha_1 L_{e_1} + \ldots + \alpha_n L_{e_n}$ is
identified with the operator matrix
\[
\left[
\begin{matrix}
0 & & & & \alpha_n T_z \\
 \alpha_1 T_z & 0 & & &  \\
 & \alpha_2 T_z & 0 & &  \\
 & & \ddots & \ddots &  \\
 & & & \alpha_{n-1}T_z & 0
\end{matrix}\right].
\]

Writing $H^\infty (z^n)$ for the subalgebra of $H^\infty$ arising from functions of the form
$h(z^n)$ with $h$ in $H^\infty$, the algebra $\fL_{C_n}$ is readily identified with the matrix
function algebra
\[
\left[
\begin{matrix}
H^\infty (z^n) & z^{n-1}H^\infty (z^n)& \hdots & z H^\infty (z^n) \\
 z H^\infty (z^n) & H^\infty (z^n) & & \vdots  \\
\vdots &   & \ddots &  \\
z^{n-1} H^\infty (z^n) & \hdots &  & H^\infty (z^n)
\end{matrix}\right].
\]

If $H^\infty$ is replaced with the disc algebra, this algebra becomes the matrix function
realization of the semicrossed product $\bbC^n \times_\beta \bbZ_+$ for the cyclic shift
automorphism of $\bbC^n$. See De Alba and Peters \cite{DeAPe} for details. It follows that
$\fL_{C_n}$ is identifiable with the $\wot$-closed semicrossed product algebra $\bbC^n
\times_\beta^\sigma \bbZ_+$.
\end{eg}

\begin{eg}\label{infinitecyclegraph}
Let $C_\infty$ be the infinite directed graph $C_\infty = (\bbZ, E)$ where the edge set
$E= \{ (n,n+1): n \in \bbZ \}$, and hence the tree components for
$C_\infty$ give the two-way
infinite graph generalization of the previous example. We have the
decomposition of Fock space $\H \equiv \H_{C_\infty}$ into a direct sum of spaces $\H_k$,
$k\in\bbZ$, where $\H_k = L_{e_k}^*L_{e_k} \H = P_k \H$ is naturally identified
with $H^2$. Pictorially, the basis elements of $\H_k$ are diagonally distributed to the
south east of vertex $x_k$. With respect to the decomposition $\H =
\sum_{k\in\bbZ} \oplus \H_k = \ell^2(\bbZ) \otimes H^2$, the partial
isometry $L_{e_k}$ has the representation $e_{k+1, k} \otimes T_z$, where
$\{ e_{ij} : i,j \in \bbZ\}$ is the standard matrix unit system for the
standard basis of $\ell^2(\bbZ)$. Since $L_{e_k}^*L_{e_k}$ has the form
$e_{k,k}\otimes I$ it follows that $\fL_{C_\infty}$ is unitarily
equivalent to the operator algebra of matrices
\[
\big( a_{ij} T_z^{j-i} \big)_{i,j = -\infty}^\infty
\]
where $(a_{ij})$ is the standard matrix of an operator in
$\T_\bbZ$, the nest algebra on $\ell^2(\bbZ)$ for the nest subspaces
$\spn \{ e_j : j\leq k\}$, $k \in \bbZ$. Once again the algebraic
structure of $\fL_{C_\infty}$ becomes clear in this representation since
$\fL_{C_\infty}$ is isomorphic to $\T_\bbZ$ by a $\wot$-$\wot$-continuous
isometric isomorphism. In particular, the Jacobson radical of
$\fL_{C_\infty}$ is determined by Ringrose's criterion \cite{nestalg} and the radical is not
closed in the weak operator topology.

On the other hand, as before, it is less clear in this representation how
to describe the commutant of $\fL_{C_\infty}$ and other spatial structure
of the algebra. Our earlier Fock space arguments show that the commutant
is unitarily equivalent to $\flgee$ where $G$ is the transpose of
$C_\infty$, and so, as with the cycle algebras $\fL_{C_n}$, the algebra
$\fL_{C_\infty}$ is isomorphic to its commutant.
%In fact the commutant of
%$\fL_{C_\infty}$ can be identified as the $\wot$-closed operator algebra
%on $\ell^2(\bbZ)\otimes H^\infty$ generated by the projections
%\[
%E_{kk} = \big( g_{k,k}\otimes e_{11}\big) \oplus \big( g_{k+1, k+1}\otimes
%e_{22} \big)\oplus \ldots,
%\]
%for $k\in\bbZ$, and the partial isometries
%\[
%E_{k+1,k} = \big( g_{k,k}\otimes e_{2,1} \big) \oplus \big( g_{k+1, k+1}
%\otimes e_{3,2} \big) \oplus \ldots .
%\]
\end{eg}

\subsection{Amalgamated Free Products}\label{sS:amal}
%%%%%%%%%%%%%%%%%%%%%%%%%%%%

We now describe how free semigroupoid algebras may be presented and decomposed in terms of
free products and amalgamated free products. The discussion here is independent of the rest of
the paper. However, the gauge automorphisms of $\flgee$ are of interest in their own right and
in the case of $\fL_n$ we use the gauge automorphisms in the proof of
Theorem~\ref{reflexive}.

Let $G$ be a countable directed graph. It is convenient now to identify the Fock space elements
$\xi_{e_{i_1}\cdots e_{i_k}}$ with tensor products $\eta_{e_{i_1}}\otimes \ldots \otimes
\eta_{e_{i_k}}$ where $\{ \eta_{e_{i_1}}, \ldots, \eta_{e_{i_n}}\} $ is a fixed orthonormal
basis for $\bbC^n$. This allows us to define certain automorphisms of $\flgee$ in a succinct
manner in terms of their action on these admissable tensors. Thus, suppose that for each vertex
pair $x_i$, $x_j$ we have a complex unitary matrix $U_{i,j}$ of size $a_{ij}\times a_{ij}$
where, as before, $a_{ij}$ is the number of edges in $G$ from $x_j$ to $x_i$. We view
$U_{i,j}$ as a unitary on
\[
\bbC^{a_{ij}} = \spn\{ \eta_e : e=x_iex_j \}.
\]
Now we may define a unitary $U = \sum\oplus U_{ij}$ on $\bbC^n$ where $U\eta_e = U_{ij}
\eta_e$ if $e = x_i ex_j$. For convenience write $U_e$ for $U_{ij}$ for any edge (there may be
several) with $e = x_i ex_j$. By linear extension there is an associated operator $\td{U}$ on the
Fock space $\H_G$ for $G$ which fixes vacuum vectors such that
\[
\td{U} \eta = \big( U_{e_1}\eta_{e_1} \big) \otimes \ldots \otimes \big( U_{e_k}\eta_{e_k}
\big)
\]
for all admissable tensors $\eta_{e_1}\otimes \ldots \otimes \eta_{e_k}$. Note, of course, that in
the expansion of the tensor $\td{U} \eta$ one obtains a linear combination of admissable tensors
and hence  $\td{U}$ is well defined as a unitary on $\H_G$.

In the free semigroup case the unitaries $\td{U}$ correspond to the so-called gauge unitaries in
the formulation of quantum mechanics. The new aspect here is that gauge unitaries are also
available for non-loop edges. Plainly the gauge unitaries $\td{U}$ respect the natural grading of
Fock space. Moreover, one can verify, by considering actions on the generators $L_e$, that the
map $\Theta_U$ on $\flgee$ defined by $\Theta_U(A) = \td{U}^* A \td{U}$ yields an injective
endomorphism of $\flgee$. Since $U$ is unitary this endomorphism is in fact an automorphism
which we refer to as a {\it gauge automorphism} of $\flgee$.

The tensor presentation of Fock space above is a simple variation of the well known presentation
of Fock space for the free semigroup algebra $\fL_n$ in the form
$
\H_n = \bbC \xi \oplus \big( \bigoplus_{N\geq 1} \K^{\otimes N} \big).
$
where $\K$ is a Hilbert space with basis $\{ \eta_{e_{i_1}}, \ldots, \eta_{e_{i_n}}\} $ indexed
by the loop edges for the single vertex graph for $\fL_n$. In general we can write the Fock space
for $G$ as
\[
\H_G = \Big( \sum_{x\in V(G)} \oplus \bbC \xi_{x}\Big)  \oplus\left(
\bigoplus_{N\geq 1} \K^{(N)} \right)
\]
where $\K^{(N)}$ is spanned by the admissable tensors of length $N$.

On the other hand, $\H_n$ can also be viewed as an $n$-fold free product of the Hilbert space
of functions $H^2$, with distinguished vector {\bf 1}. In fact, in general, if $(\H_i, \xi_i)$ is a
family of
Hilbert spaces with distinguished unit vectors, then their Hilbert space free product $(\H,\xi)$,
denoted $\ast_{i\in \I}(\H_i,\xi_i)$, has the form
\[
\H = \bbC \xi \oplus \bigoplus_{N\geq 1} \left( \bigoplus_{i_1\neq i_2 \neq \ldots \neq i_N}
\stackrel{\circ}{\H_{i_1}} \otimes \ldots \otimes \stackrel{\circ}{\H_{i_N}}\right)
\]
where $\stackrel{\circ}{\H_{i}} = \H \oplus \bbC \xi_i$. Such a free product Hilbert space
allows for the construction of a faithful representation
$
\pi = \ast_{i\in \I}  \pi_i  :  \ast_{i\in \I}\, \A_i \rightarrow \bofh
$
of an (algebraic) free product of unital operator algebras $\A_i$, each represented on a Hilbert
space $\H_i$ by a faithful representation $\pi_i$. Viewing the operator algebras $\A_i$ as
represented on $\H_i$ we suppress the notation $\pi_i$ and define $\ast_{i\in \I} \,\A_i$ in the
category of $\wot$-closed operator algebras to be the $\wot$-closure of the algebraic free
product in the representation $\pi$. For more details see \cite{DNV}. In this way we can define
the $\wot$-closed algebra $\ast_{i\in \{ 1, \ldots, n\}}H^\infty$, where $H^\infty$ is represented
on $H^2$ with distinguished vector {\bf 1}, and this can be shown to coincide with $\fL_n$.
However, in the special case of free products of the operator algebras $\fL_n$ one may take a
slightly more explicit approach: Let $\K = \K_1 \oplus \K_2 = \bbC^{n_1} \oplus \bbC^{n_2}$
be a decomposition associated with a partition of the basis $\{ \eta_{e_{i_1}}, \ldots,
\eta_{e_{i_n}}\} $, where $n = n_1 + n_2$. Then $\fL_n$ acts on the unrestricted Fock space
%\begin{eqnarray*}
%\H_N &=& \bbC \xi \oplus \left( \bigoplus_{n\geq 1} \K^{\otimes n} \right)
%= \bbC \xi \oplus \left( \bigoplus_{n \geq 1} (\K_1 \oplus\K_2)^{\otimes n} \right) \\
%&=& \bbC \xi \oplus \bigoplus_{n\geq 1} \left( \bigoplus_{i_1\neq i_2 \neq \ldots \neq i_n}
%\Big(
%\big[ \bigoplus_{p_1\geq 1} \K_{i_1}^{\otimes p_1} \big] \otimes \ldots
% \otimes \big[ \bigoplus_{p_n\geq 1} \K_{i_n}^{\otimes p_n} \big]
%\Big)\right).
%\end{eqnarray*}
$\H_n$, which can be recognized as the free product $\big( \H_{n_1},\xi \big)\ast \big(
\H_{n_2},
\xi\big)$ defined above. (This is a special case of (\ref{freeprodspace}) below.) Plainly
$\fL_{n_1}$ and $\fL_{n_2}$ are
naturally represented on this space by viewing each $L_e$ as an appropriate creation operator.
One can check that this representation agrees with our general free product representation $\pi$,
and it follows that $\fL_{n_1}\ast\fL_{n_2}$ is unitarily equivalent to $\fL_n$.

Suppose now that $G_1$ and $G_2$ are finite directed graphs with a single identified vertex
$x$ say, and let $G = G_1 \sqcup_x G_2$ be the amalgamated graph. We want to identify
$\flgee$ with a $\wot$-closed amalgamated free product $\fL_{G_1} \ast_{\bbC P_x}
\fL_{G_2}$, where $\bbC P_x$ is the common subalgebra. Rather than placing this construction
in a general context, let us take advantage of the explicit nature of the algebras $\flgee$, as
defined by creation operators, to give a direct definition of $\fL_{G_1} \ast_{\bbC P_x}
\fL_{G_2}$. This definition, as before, relies on specifying an appropriate Hilbert space on
which to represent $\fL_{G_1}$ and $\fL_{G_2}$. The construction is similar to the free
semigroup case except that only admissable tensors corresponding to finite paths are considered.
Thus we have
\[
\H_{G_1} = \spn \{ \xi_{x_i}\} \oplus \bigoplus_{N\geq 1}\K_1^{(N)} \qand
\H_{G_2} = \spn \{ \xi_{y_j}\} \oplus \bigoplus_{N\geq 1}\K_2^{(N)}.
\]
With $x_1 = y_1 = x$ understood we write $\K_1^{(N)} \otimes_x \K_2^{(M)}$ to denote the
Hilbert space spanned by basis elements
\[
\big( \eta_{e_{i_1}} \otimes \ldots \otimes \eta_{e_{i_N}}\big) \otimes
\big( \eta_{f_{j_1}} \otimes \ldots \otimes \eta_{f_{j_M}}\big)
\]
where $e_{i_N} = e_{i_N}x$ and $f_{j_1} = x f_{j_1}$. Now we may define
\begin{eqnarray}\label{freeprodspace}
\H &=& \H_{G_1} \ast_x \H_{G_2} \equiv \spn \{ \xi_{x_i}, \xi_{y_j} \} \oplus
\bigoplus_{N\geq 1}
\end{eqnarray}
\begin{eqnarray*}
& & \left( \bigoplus_{i_1\neq i_2 \neq \ldots \neq i_N} \Big( \big[ \oplus_{p_1 \geq 1}
\K_{i_1}^{(p_1)} \big] \otimes_x \ldots \otimes_x \big[ \oplus_{p_N \geq 1}
\K_{i_N}^{(p_N)}
\big] \Big) \right).
\end{eqnarray*}
There are natural representations $\pi_i : \fL_{G_i} \rightarrow \bofh$, which agree on the
common projection $P_x$ and we define $\fL_{G_1}\ast_{\bbC P_x} \fL_{G_2}$ to be the
generated $\wot$-closed operator algebra. By  simple manipulations as in the free semigroup
case, we can identify this operator algebra with $\flgee$.
Thus we have the following.

\begin{thm}\label{freeprodthm}
Let $G_1$ and $G_2$ be finite directed graphs with a single identified vertex $x$ and let $G=
G_1 \sqcup_x G_2$ be the amalgamated graph. Then $\flgee$ is naturally unitarily equivalent to
the amalgamated free product algebra $\fL_{G_1}\ast_{\bbC P_x} \fL_{G_2}$.
\end{thm}

One can also view  $\fL_{G_1}\ast_{\bbC P_x} \fL_{G_2}$ as arising from the amalgamation
over $x$ of the left regular representations $\lambda_{G_1}$, $\lambda_{G_2}$ of the free
semigroupoids $\bbF^+(G_1)$, $\bbF^+(G_2)$. That is, from $\lambda_{G_1}\ast_x
\lambda_{G_2}$, appropriately defined, where $x$ is an identified unit of $G_1$ and $G_2$.
As in the group case (see \cite{DNV}) $\lambda_{G_1}\ast_x \lambda_{G_2}$ identifies
naturally with $\lambda_{G_1\sqcup_x G_2} = \lambda_G$.
We can now revisit the $\fL_{G^\prime}$ of
Example~\ref{easymtxfnalg} and see that it is unitarily equivalent to $H^\infty \ast_{\bbC P_x}
\fL_{C_2}$. More generally, if $G$ is the graph formed by joining $C_n$ and $C_m$ at a
single vertex $x$ then
\[
\flgee \simeq \big( \bbC^n \times^\sigma_\beta \bbZ_+ \big) \ast_{\bbC P_x}
\big( \bbC^m \times^\sigma_\beta \bbZ_+ \big).
\]

%%%%%%%%%%%%%%%%%%%%%%%%%%%%%%%%%%%%%%%%%%%%%%
%%%%%%%%%%%%%%%
\section{Reflexivity and Eigenvectors}\label{S:invariant}
%%%%%%%%%%%%%%%%%%%%%%%%%%%%%%%%%%%%%%%%%%%%%%
%%%%%%%%%%%%%%%

We first give an elementary proof that the algebras $\flgee$ are
reflexive. Recall that given an operator algebra $\fA$ and a
collection of subspaces $\L$, the subspace lattice $\Lat \fA$
consists of those subspaces left invariant by every member of
$\fA$, and the algebra $\Alg \L$ consists of all operators which
leave every subspace of $\L$ invariant. Every algebra satisfies
$\fA \subseteq \Alg \Lat \fA$ and an algebra is {\it reflexive} if
$\fA = \Alg \Lat \fA$.

\begin{thm}\label{reflexive}
$\flgee$ is reflexive.
\end{thm}

\Prf Let $A \in \Alg \Lat (\flgee)$. Then $A( \ol{R\H_G})
\subseteq \ol{R\H_G}$ for all $R\in \frgee$ since $\flgee =
\frgee^\prime$. Suppose that $A = AP_x$ for some $x\in V(G)$. For
$u=yu\in \fngee$, the subspace $R_u \H_G$ is the set of vectors in
$\H_G$ of the form $\sum_{s(w)=y} \beta_w \xi_{wu}$. Hence $A
\xi_u = A R_u \xi_y = A R_u P_x \xi_y =0$ when $x\neq y$, and
otherwise $A\xi_u$ has expansion $\sum_{s(w)=x} \alpha_w^u
\xi_{wu}$ when $r(u)=x$. We claim that the scalars $\alpha_w^u$
are independent of $u$; that is, if $r(u)=x=r(v)$, then
$\alpha_w^u =\alpha_w^v$ for all $w=wx$.

Assume the claim holds for the moment. Then there are scalars $\{
\alpha_w \}_{w\in \fngee}$, where $\alpha_w = 0$ for $s(w) \neq
x$, with $A\xi_u = \sum_{s(w)=x } \alpha_w \xi_{wu}$ for all
$u=xu\in \fngee$ and $A \xi_u =0$ otherwise. For such an operator,
it is easy to check that $A$ belongs to $\frgee^\prime = \flgee$.
In the general case, given $A \in \Alg \Lat (\flgee)$ we may write
$A = \sum_{x\in V(G)} AP_x$,  the sum converging in the strong
operator topology for the infinite vertex case. Then the above
argument can be applied to show each $AP_x$ belongs to
$\frgee^\prime = \flgee$, and hence $A$ is in $\flgee$ as
required. Thus, we will prove the theorem by verifying the claim.

There are two cases to consider. For the first case let us suppose
that there is an edge $v$ with $v=xvy$ and $y\neq x$. Then the
range $\M$ of $R_x + R_v$ is spanned by the set of vectors
$\{\xi_{wx} + \xi_{wv} : s(w)=x \}$ and the vectors in this set
are pairwise orthogonal. Since $\M \in \Lat \flgee$ it follows
that
\[
A(\xi_x + \xi_v) = A(R_x + R_v) \xi_x = \sum_{s(w)=x} h_w
(\xi_{wx} + \xi_{wv})
\]
for some choice of scalars $h_w$. But $A\xi_x$ and $A\xi_v$ are
given in terms of the coefficients $\alpha_w^x$, $\alpha_w^v$
respectively and we conclude that $\alpha_w^x = h_w = \alpha_w^v$
for all $w$.

For the remaining case we have $u=xu$ only if $u=xux$ and we see
that the compression algebra $P_x \flgee|_{P_x\H_G}$ is unitarily
equivalent to $\fL_H$ where $H$ is a single vertex subgraph of
$G$. Thus $\fL_H$ is unitarily equivalent to $\fL_n$ for some
$n\geq 1$ or $\fL_H$ is unitarily equivalent to $\bbC$, and these
algebras are known to be reflexive. We may assume then that $H\neq
G$ and also that $\alpha_w^v = \alpha_w^x$ for all words $v,w$
with $r(v)=x = r(w)$. Suppose that $w' = y w' x$ with $y\neq x$.
It remains to show that $\alpha_{w'}^v = \alpha_{w'}^x$.

Suppose first that $w'$ is not of the form $w_1h$ where $h$ is a
path in $H$ with $|h|\geq 1$. Consider the restriction operator
$A_{w'} = L_{w'}^* A|_{\H_H}$. We show that $A_{w'}$ is in
$\fL_H$. To this end, let $\M\in \Lat \fL_H$ and define $\td{\M} =
\bigvee_{s(w)=x} L_w \M$ in $\Lat \flgee$. Then
\[
A_{w'} \M = L_{w'}^* A \M \subseteq L_{w'}^* \td{\M}.
\]
In view of our assumption on $w'$ we have that $L_{w'}^* L_w$ is
non-zero only if $w=w'h$ with $h$ a path in $H$. Thus
\[
L_{w'}^* \td{\M} \subseteq \bigvee_{w=xwx} L_w \M \subseteq \M.
\]
We have shown that $A_{w'}\in \Alg \Lat \fL_H$ and so $A_{w'}\in
\fL_H$. Hence there are scalars $\alpha_h$ such that
\[
A_{w'} \xi_v = \Big( \sum_{h=xhx} \alpha_h L_h \Big)\, \xi_v
\]
for all paths $v$ in $H$. Thus we have
\begin{eqnarray*}
\sum_h \alpha_h \xi_{w'hv} &=& L_{w'} \big( A_{w'} \xi_v \big) \\
&=& L_{w'} L_{w'}^* A \xi_v \\
&=& L_{w'} L_{w'}^* \Big( \sum_w \alpha_w^v \xi_{wv} \Big) \\
&=& \sum_h \alpha_{w'h}^v \xi_{w'hv}.
\end{eqnarray*}
This shows that $\alpha_{w'h}^v = \alpha_h$ for all $h=xhx$. In
particular, we obtain $\alpha_w^v = \alpha_w^x$ for all words $w$
with $w=ywx$, $y\neq x$, as desired. \bx

We next compute the eigenvalues for $\flgee^*$, by which we mean, with modest
abuse of terminology, the values
$\lambda = (\lambda_e)_{e\in\edges(G)}$ in $\bbC^n$, where $n$ is the cardinality of $E(G)$,
for which there is a
unit vector $\xi$ in $\H$ such that $L_e^* \xi = \ol{\lambda}_e \xi$ for all $e
\in\edges(G)$.  Since the $L_e$ are partial isometries with pairwise orthogonal
final projections, we have
\[
\sum_{e\in\edges(G)} |\lambda_e|^2 = \sum_{e\in\edges(G)}
(L_eL_e^*\xi,\xi) \leq ||\xi||^2 = 1.
\]
In the free semigroup case, the open unit $n$-ball
\[
\bbB_n = \Big\{ \lambda =  (\lambda_e)_{e\in\edges(G)}
\,\, : \,\, \sum_{e\in\edges(G)} |\lambda_e|^2 < 1 \Big\}
\]
forms the set of eigenvalues for $\fL_n^*$ \cite{AP,DP1}.
In general the eigenvalues for $\flgee^*$ form a proper
subset of the unit $n$-ball, with structure determined by  lower dimensional unit balls.

The set of eigenvalues for $\flgee^*$ will be described explicitly in terms
of  $G$. We begin by pointing out a special case which is quite
different from the semigroup case.

\begin{prop}\label{noevalues}
If $G$ has no loop edges,
$a_{xx}=0$ for all $x$ in $V(G)$, then $\lambda =
\vec{0} \in\bbC^n$ is the only eigenvalue for $\flgee^*$.
\end{prop}

\Prf
In this case $L_e^2 = L_{e^2} =0$ for
all $e\in\edges(G)$, as no words of the form $e^2w$ belong to $\fngee$. As every eigenvalue
$\lambda = (\lambda_e)_{e\in\edges(G)}\in\bbC^n$ for
$\flgee^*$ is determined by  equations $L_e^* \xi = \ol{\lambda}_e
\xi$, it follows that $\lambda = \vec{0}$.
\bx

We require some extra notation to state the following theorem. We shall say  the `tree
top graph' associated with a directed graph $G$ (the top two levels of the tree components for
$\H_G$) has a {\it standard ordering} if  the
{\it saturation} $G_x$ (the set of all paths which start at $x$) at every vertex $x$ in $V(G)$  has
all its  edges which finish at $x$ lying to the left of all other edges in $G_x$. Also, we shall
use the
notation $0_k = \vec{0}\in\bbC^k$. Further, we will assume the vertices of $G$ are given by
$x_1, \ldots , x_M$.  Recall that $A = (a_{x_ix_j})\equiv  (a_{ij})$ is the transition matrix
associated
with $G$, where  $a_{ij}$ is the number of directed edges from vertex $x_j$ to vertex $x_i$.
Finally, write $W_i$ for the set of all words in edges which are loops at
vertex $x_i$, and put $P_i = P_{x_i}$, $Q_i = Q_{x_i}$.

\begin{thm}\label{eigenvalues}
Let $G$ be a countable directed graph with tree top graph having a standard edge ordering.   Let
$A = (a_{ij})$ be the
transition matrix for the graph, and let $k_i = \sum_{j=1}^M a_{ij}$ for
$1 \leq i \leq
M$. Then
\begin{itemize}
\item[$(i)$] Every eigenvector for $\flgee^*$ belongs to $P_i\H_G$ for some
$1 \leq i \leq M$.
\item[$(ii)$] The unit eigenvectors supported on $P_i\H_G$ are scalar multiples of the
vectors
\[
\mulam = (1 - ||\lambda||_2^2)^{1/ 2} \sum_{w\in W_i} \ol{w(\lambda)}\xi_w
\]
where $\lambda = (\lambda_e)_{e\in\edges(G)}\in\bbC^n$, belongs to
the set
\[
\bbB_{G_{x_i}} \equiv (0_{k_1},\ldots, 0_{k_{i-1}}, \bbD_i, 0_{k_{i+1}}, \ldots,
0_{k_M})\in\bbC^n,
\]
where $\bbD_i \equiv (\bbB_{a_{ii}}, 0_{k_i - a_{ii}})\in\bbC^{k_i}$. Also,
\[
\nu_{\lambda,i} = (1- || \lambda ||_2^2)^{1/ 2} \Big( I -
\sum_{e =x_i e x_i} \ol{\lambda}_e L_e \Big)^{-1} \xi_{x_i}.
\]
\item[$(iii)$] The eigenvectors $\mulam$ are supported on
$Q_iP_i\H_G$ and are $Q_i\H_G$-cyclic for $\flgee$.
They satisfy
\[
L_e^* \mulam = \ol{\lambda}_e \mulam
\]
and if $L_G$ is the $n$-tuple $L_G= (L_e)_{e\in\edges(G)}$, then
$(p(L_G)\mulam, \mulam) = p(\lambda)$ for every polynomial $p = \sum_w a_w
w$ in the semigroupoid algebra $\bbC \,\fngee$. This extends to the map
$\varphi_{\lambda,i} (A) = (A \mulam, \mulam)$, which is a
$\wot$-continuous multiplicative linear functional on $\flgee$.
\end{itemize}
\end{thm}

\Prf
Let  $\nu$ be an eigenvector for $\flgee^*$.  Then $\nu$ is an eigenvector for
the projections $P_1, \ldots, P_M$, but the only eigenvalues for a
projection are $0$ or $1$. Thus, as the $P_i$ have pairwise orthogonal
ranges summing to the identity, there is a unique $i$ with $\nu = P_i
\nu \in P_i \H_G$.

There are scalars $\lambda_e$ such that $L_e^* \nu = \ol{\lambda_e} \nu$ for all $e$. If $\nu =
\sum_w \ol{a_w} \xi_w$, then
\[
\sum_w \ol{\lambda_e} \ol{a_w} \xi_w = \ol{\lambda_e} \nu = L_e^* \nu =
\sum_{w =ev} \ol{a_{ev}} \xi_v
\]
and so $\ol{\lambda_e}\ol{a_w} = \ol{a_{ew}}$, and thus $\ol{a_w} =
\ol{w(\lambda)}$.
However, there will typically be $\lambda_e$ equal to  zero. If $e$ is an edge with distinct
initial and final vertices, then
$L^2_e = 0$ because there are no words in $\fngee$ of the form
$e^2w$. Hence for such edges we have $\lambda_e = 0$. Further, let $e$
be an edge in $G$ with initial  vertex distinct from $x_i$; that is, $e=ex_j$ for
$j\neq i$. Then $L_e = L_e P_j$, and
\[
P_i (\ol{\lambda_e} \nu) = \ol{\lambda_e} \nu = L_e^* \nu = P_j L_e^* \nu
= P_j (\ol{\lambda_e} \nu),
\]
which can only happen if $\lambda_e = 0$. Therefore we have shown that
the eigenvalues corresponding to eigenvectors supported on $P_i \H_G$ must
belong to the set $\bbB_{G_i}$ (As $\nu$ has finite norm, it follows that $||\lambda|| < 1$), and
the vectors are supported on
$\spn \{ \xi_w : w\in W_i \}$.

Now given $\lambda = (\lambda_e)_{e\in\edges(G)}$ in $\bbB_{G_i}$,
%\begin{eqnarray*}
%\sum_{w=x_iw} |w(\lambda)|^2 &=& \sum_{s\geq 0} \sum_{v=j_1\cdots j_s,
%\; 1\leq j_l \leq a_{ii}} |v(\lambda_1, \ldots, \lambda_{a_{ii}})|^2 \\
%&=& \sum_{s\geq 0} \big( \sum_{j=1}^{a_{ii}} |\lambda_j|^2 \big)^s
%=(1- ||\lambda||^2)^{-1} < \infty.
%\end{eqnarray*}
we have
\[
\Big|\Big| \sum_{e\in W_i} \ol{\lambda_e} L_e\Big|\Big|^2 = \Big|\Big| \sum_{e\in W_i}
|\lambda_e|^2 L_e^* L_e \Big|\Big| \leq \sum_e |\lambda_e|^2 = ||\lambda||^2 < 1,
\]
so that $I - \sum_e \ol{\lambda_e} L_e$ is invertible, and its inverse is
given by the power series
\[
\left( I - \sum_e \ol{\lambda_e} L_e \right)^{-1} = \sum_{k \geq
0} \left( \sum_e \ol{\lambda_e} L_e \right)^k = \sum_{r(w)=x_i}
\ol{w(\lambda)} L_w.
\]
Thus the second identity for $\mulam$ follows, and from this it is clear
that $\mulam$ is $Q_i\H_G$-cyclic for $\flgee$.

The vectors $\mu_{\lambda,i} = (1- ||\lambda||^2)^{-1 / 2} \mulam$ satisfy
\[
L_e^* \mu_{\lambda,i} = L_e^* \sum_{r(w)=x_i} \ol{(ew)(\lambda)}
\xi_{ew} = \ol{\lambda_e} \sum_{r(w)=x_i} \ol{w(\lambda)} \xi_w =
\ol{\lambda_e} \mu_{\lambda,i},
\]
and we also have $(L_w \nulam, \nulam) = w(\lambda) ||\nulam||^2 =
w(\lambda)$, which easily extends to polynomials by linearity. It is clear
that $\varphi_{\lambda,i}$ is multiplicative and $\wot$-continuous.
\bx

\begin{rem}\label{evectorremark}
There is an analogous version of this result for the eigenvectors of $\frgee^*$, with the operators
$\{R_e, Q_x \}$ in place of  $\{ L_e, P_x \}$. This fact is used in the proof of
Theorem~\ref{autothm}. Observe that the form of non-vacuum eigenvectors $\nu_{\lambda, x}
= P_x Q_x \nu_{\lambda, x}$ for $\frgee^*$ are also determined by the $a_{xx}$ loop edges
over vertex $x$.
\end{rem}

%%%%%%%%%%%%%%%%%%%%%%%%%%%%%%%%%%%%%%
\section{Beurling Theorem and Partial Isometries}\label{S:beurling}
%%%%%%%%%%%%%%%%%%%%%%%%%%%%%%%%%%%%%%

In this section we establish a Beurling-type invariant subspace theorem for $\flgee$.
As a consequence we obtain a structure theorem for partial isometries in $\flgee$, and an
inner-outer factorization for elements of $\flgee$.

We will say that a non-zero subspace $\W$ of $\H$ is {\it
wandering} for $\flgee$ if the subspaces $L_w\W$ are pairwise
orthogonal for distinct $w$ in $\fngee$. Observe that  the partial
isometries $L_w$ with $w\in\fngee$ include the vertex projections
$P_x = L_{x}$. Further, since the $L_e$ are partial isometries we
cannot `peel off' $L_e$'s when comparing the subspaces $L_w\W$, as
is done in the case of isometries with orthogonal ranges.
Nonetheless, equations such as $L_w^* L_w = P_{s(w)}$ give us a
computational device for this comparison process.

Every wandering subspace $\W$ generates an
$\flgee$-invariant subspace
\[
\flgee [\W] = \sum_{w\in\fngee} \oplus L_w \W.
\]
Every $\flgee$-wandering vector $\zeta$ generates the cyclic
invariant subspace $\flgee[\zeta]$. We will say $\flgee[\zeta]$ is
a {\it minimal cyclic subspace} if $P_x \zeta = \zeta$ for some
$x\in V(G)$. It is easy to see that  given a wandering vector
$\zeta$, each vector $P_x \zeta$ which is non-zero is wandering as
well.

\begin{thm}\label{beurlingthm}
Every invariant subspace of $\flgee$ is generated by a wandering subspace,
and is the direct sum of minimal cyclic subspaces generated by wandering
vectors. Every minimal cyclic invariant subspace generated by a
wandering vector is the range of a partial isometry in $\frgee$, and the choice
of partial isometry is unique up to a scalar multiple.
\end{thm}

\Prf
Let $\M$ be a non-zero invariant subspace for $\flgee$ and form the
subspace
\[
\W = \M \ominus \left( \sum_{e\in\edges(G)} \oplus L_e \M \right).
\]
First note that $\W$ is a wandering subspace for $\flgee$. To see
this, let $\xi$ and $\eta$ be vectors in $\W$ and let $v$, $w$
belong to $\fngee$. Consider the inner product $(L_w\xi, L_v
\eta)$. This is clearly zero if $v$ and $w$ are distinct units in
$V(G)$. If $w = x$ is a unit and $|v| \geq 1$, then $(P_x \xi, L_v
\eta)=0$ if $r(v)\neq x$, and otherwise $(\xi,L_v\eta)=0$ by the
definition of $\W$. If $v$ and $w$ are non-units with differing
left most letters, then $(L_w\xi, L_v\eta)=0$. Otherwise we would
have $w=ew_1$ and $v=ev_1$ so that $ (L_w\xi, L_v \eta)=
((L_e^*L_e)L_{w_1}\xi, L_{v_1}\eta). $ Since $L_e^*L_e=P_x$ for
some $x$, we may repeat this argument to  show  this inner product
is always zero.

We claim that $\M= \flgee[\W]$. Let $\N$ be the orthogonal complement of
$\flgee[\W]$ inside $\M$. Let $\eta\in\N$ and let $\xi=\sum_w
L_w \zeta_w$, with each $\zeta_w\in\W$, belong to $\flgee[\W]$. Then for
all  $e=ex$ in $E(G)$ we have
\[
(L_e\eta, \xi) = \sum_w (L_e\eta, L_w \zeta_w)
= \sum_{u\in\fngee}(\eta, P_xL_u \zeta_w)
= 0,
\]
from the definition of $\N$. Further, $(P_x \eta, \xi) = (\eta,
P_x \xi) = 0$ since $P_x \xi \in \flgee[\W]$. Thus it follows that
$\N$ is invariant for $\flgee$. Now let $\eta$ belong to the
orthogonal complement of $\sum_e \oplus L_e \N$ inside $\N$. As
$\eta$ belongs to $\N$, we know that $(\eta, L_e\xi)=0$ for all
$\xi\in\flgee[\W] \subseteq \M$. It follows that $\eta$ belongs to
$\W$. Indeed, let $\zeta\in\M$ and put $\zeta = \zeta_1 + \zeta_2$
with $\zeta_1\in\N$ and $\zeta_2 \in \flgee[\W]$. Then
$L_e\zeta_1\in\N$ and $L_e\zeta_2 \in \flgee[\W]$ so that $ (\eta,
L_e\zeta) = (\eta, L_e \zeta_1) + (\eta, L_e\zeta_2) =0. $ Thus we
have established that $\eta$ is a vector in $\N$ which is also in
$\W \subseteq (\N)^\perp$, whence $\eta = 0$. In particular, we
have $\N = \sum_e \oplus L_e \N$. Finally, assume that $\N \neq \{
0 \}$ and let $k_0$ be minimal with $E_{k_0} \N \neq \{0\}$. Then
\[
E_{k_0}\N \subseteq \sum_e E_{k_0} L_e \N = \sum_e L_e E_{k_0-1}\N =0.
\]
This contradiction yields $\N = \{0 \}$, and hence $\M = \flgee[\W]$ as
claimed.

Next we observe that $\M$ is the direct sum of cyclic subspaces. First
note that $P_x\W \subseteq \W$ for all $x\in V(G)$ as noted in the discussion preceding the
theorem. For
each $x$, let $\{ \zeta_{x,k_x}\}$ be an orthonormal basis for $P_x \W
\subseteq \W = \sum_{x\in V(G)} \oplus P_x\W$. Then
\[
\flgee[\W] = \sum_x \sum_{k_x} \oplus \flgee[\zeta_{x,k_x}].
\]
The right-hand side of this identity is contained in $\M =
\flgee[\W]$ by $\flgee$-invariance. On the other hand, note that
for distinct wandering vectors $\zeta_{x,k_x}$ and $\zeta_{y,l_y}$
in $\W$ (with $x$ and $y$ not necessarily distinct), the cyclic
subspaces $\flgee[\zeta_{x,k_x}]$ and $\flgee[\zeta_{y,l_y}]$ are
perpendicular. Lastly, it is clear that vectors in $\flgee[\W]$
belong to the sum on the right side.

The subspace $\M$ is minimal cyclic if and only if $\W$ is
one-dimensional and there is an $x$ with $P_x \W = \W$ and $P_y\W
= 0$ for $y\neq x$. Consider $\M = \flgee[\zeta]$, where $P_x
\zeta = \zeta$ is a unit $\flgee$-wandering vector. Define a
linear transformation $R_\zeta$ on $\H_G$ by the rule
$R_\zeta\xi_w = L_w \zeta$ for $w\in\fngee$. Then for $w\neq v$ in
$\fngee$ we have $(R_\zeta \xi_w, R_\zeta \xi_v)= (L_w\zeta, L_v
\zeta)=0$ because $\zeta$ is wandering. Further, when $L_w^*L_w
\neq P_x$ a similar computation shows that $||R_\zeta \xi_w||^2 =
0$. Moreover, if $L_w^*L_w = P_x$ we see that $||R_\zeta \xi_w||^2
= ||P_x \zeta||^2=1$. Thus it follows that the operator $R_\zeta$
is a partial isometry with range equal to $\M$ by design. Finally,
observe that for each edge $e$ and $w$ in $\fngee$
\[
R_\zeta L_e \xi_w = R_\zeta \xi_{ew} = L_{ew} \zeta = L_e R_\zeta\xi_w.
\]
Similarly, $R_\zeta P_y \xi_w =L_{yw} \zeta = P_y R_\zeta \xi_w$
and we have $R_\zeta\in\flgee^\prime = \frgee$ as required.

To verify the uniqueness assertion, suppose $P_x \zeta = \zeta$ is
a unit $\flgee$-wandering vector and   $\M = \flgee [\zeta] =
R_\zeta \H_G =  R\H_G$ is the range of another partial isometry
$R$ in $\frgee$. We claim that $R = \lambda R_\zeta$ for some
$|\lambda|=1$. First observe that the  vectors  $R\xi_w = L_w
R\xi_x$, where $s(w)=x$, form an orthonormal basis for $\M=
R\H_G$. This follows from Corollary~\ref{stdformcor} since the
initial projection $R^*R =Q_x$. (Note that
Corollary~\ref{stdformcor} relies on Theorem~\ref{stdformthm},
which in turn uses part of the proof of this theorem, but that
there is no circular logic. The proof of Theorem~\ref{stdformthm}
does not use the uniqueness from this theorem.) In particular,
clearly $R\xi_x$ belongs to the wandering subspace $\W = P_x \W =
\spn \{ \zeta \}$, and hence $R\xi_x = \lambda \zeta = \lambda
R_\zeta \xi_x$ for some $|\lambda| = 1$. Thus for $w=wx$ we have
$R\xi_w = L_w R\xi_x = \lambda R_\zeta \xi_w$, and it follows that
$R = RQ_x = \lambda R_\zeta Q_x = \lambda R_\zeta$. \bx

\begin{rem}\label{beurlingremark}
There is obviously an analogue of this result for the invariant subspaces of
$\frgee$, where the notion of wandering is determined by the $R_e$ and $Q_x$. This is used in
Theorem~\ref{autothm} and  Section~10. We also note that
Theorem~\ref{beurlingthm} parallels the Beurling Theorem from
\cite{MS_CJM}, and gives a slight improvement for these algebras. Indeed, we have
identified
the minimal cyclic subspaces as ranges of partial isometries in the commutant algebra $\frgee$,
and shown that the decomposition in terms of minimal invariant subspaces is
unique. Further, in
Theorem~\ref{stdformthm} we  prove that all such operators have a standard  form.
\end{rem}

The range of every partial isometry $R$ in $\frgee$ is cyclic
since $R\H_G = \ol{R\flgee\xi_\phi}= \ol{\flgee R\xi_\phi}$, where
$\xi_\phi = \sum_{x_k\in V(G)} \frac{1}{k} \xi_{x_k}$. However, we
observe through the next example that these subspaces are not
necessarily minimal cyclic. This is different from the free
semigroup case \cite{DP1,Pop_beur}, where ranges of isometries are
minimal cyclic subspaces. The basic difference here is that
partial isometries in the commutant can `cross-over' between
distinct tree components.

\begin{eg}
Let $G$ be the directed graph with transition matrix $A
=\begin{sbmatrix} 1 &1 \\ 1& 0 \end{sbmatrix}$. Let $V(G) = \{x_1,
x_2 \}$ and let $e_1 =x_1 e_1 x_1$, $e_2 = x_2e_2x_1$ and $e_3=
x_1e_3x_2$. Let $R\in \frgee$ be the isometry defined by $R=
R_{e_1} + R_{e_2}$. It is an isometry because $R_{e_1} Q_2 = 0 =
R_{e_2} Q_1$, whereas the ranges of $R_{e_1}Q_1 = R_{e_1}$ and
$R_{e_2}Q_2 = R_{e_2}$ are orthogonal. The range of $R$ is a
cyclic subspace given by
\[
\M=R\H_G = \spn\{ \xi_{we_1} + \xi_{we_2} : w\in\fngee\} = \flgee[\xi_{e_1}+
\xi_{e_2}].
\]
Thus,  $L_e\M = \spn\{ \xi_{ewe_1} + \xi_{ewe_2} : w\in
\fngee\}$ for $e\in\edges(G)$, and hence
\[
\W = \M\ominus \Big( \sum_e\oplus L_e \M\Big)  = \spn\{\xi_{e_1}, \xi_{e_2}\}.
\]
Thus $\W$ is two-dimensional here with $P_i \W$ spanned by $\xi_{e_i}$ for
$i=1,2$. In particular, from the proof of Theorem~\ref{beurlingthm}, we
see that $\M$ is the direct sum of two minimal cyclic
subspaces $\M = \flgee[\xi_{e_1}] \oplus \flgee[\xi_{e_2}]$.
\end{eg}

We next derive an explicit  characterization of  partial
isometries in $\flgee$.  We begin with a computational lemma.

%In this section, by focusing on the algebras $\flgee$, we initiate an
%investigation into deciding
%when an operator algebra contains a copy of the free semigroup algebra $\fL_2$.
%Philosophically
%this is related to the $\ca$-algebra notion of deciding when a $\ca$-algebra
%contains a copy of
%the Cuntz algebra $\O_2$. More specifically, we are interested in determining
%when a
%$\wot$-closed operator algebra $\fA$ contains an algebraically free component.
%Hence we say
%that $\fA$ is {\it partially algebraically free} if there is a $\wot - \wot$
%continuous isometric
%algebraic injection of $\fL_2$ into $\fA$. If the map is also unital, we say
%$\fA$ is {\it strongly
%partially algebraically free}.
%
%For the algebras $\flgee$ in Theorems~\ref{doublecyclethm} and
%\ref{strongdoublecyclethm}
%show that these two notions can be described strictly in terms of the graph $G$.
%In establishing these results we obtain a tight characterization of the partial
%isometries in these algebras.
%We begin by pointing out a helpful computational lemma motivated by the free
%semigroup case
%\cite{DP1}.

\begin{lem}\label{stdformlemma}
For $A\in\flgee$, we have
\begin{eqnarray}\label{computlemma}
R_f^* A - A R_f^*= \big(R_f^*A\xi_x \big) \otimes \xi_x^* \qwhere
s(f) =x.
\end{eqnarray}
\end{lem}

\Prf
Let $A\in\flgee$. Then
\[
A = AI = A \left( \sum_{x\in V(G)} Q_x \right) = A \left( \sum_x
\Big( \xi_x \otimes \xi_x^* + \sum_{s(e)=x} R_e R_e^* \Big)
\right).
\]
Hence for $f\in\edges(G)$,
\begin{eqnarray*}
R_f^*A &=& \sum_x \big(R_f^*A \xi_x\big) \otimes \xi_x^* + \sum_x
\sum_{s(e)=x} R_f^*A R_e
R_e^*, \\
&=& \sum_x \big(R_f^*A \xi_x\big) \otimes \xi_x^* + AR_f^*,
\end{eqnarray*}
since $R_f^*A R_e R_e^* = R_f^* R_e A R_e^* = \delta_{f,e} A
R_f^*$. Further, $A\xi_x = Q_x A\xi_x$ so that $R_f^*A \xi_x = 0$
when $s(f) \neq x$. The identity (\ref{computlemma}) follows. \bx

\begin{thm}\label{stdformthm}
Let $V$ be a partial isometry in $\flgee$. Then
\[
V = \sum_{i\in \I} \oplus L_{\eta_i}
\]
where $\{ \eta_i \}_{i\in\I}$ are unit wandering vectors for
$\frgee$ supported on distinct $Q_x \H_G$, and where the series
converges in the strong operator topology if $\I$ is an infinite
set.
\end{thm}

\Prf Let $\eta_1, \eta_2, \ldots $ be an orthonormal basis for the
$\frgee$-wandering subspace $\W$ of $\M = V\H_G$, where each
$\eta_i$ belongs to some $Q_x \H_G$, as in the proof of the
$\frgee$ version of Theorem~\ref{beurlingthm}. (In fact, it
follows that the set $\{ \eta_i \}_{i\in\I}$ is equal to the set
of non-zero vectors amongst $\{ V\xi_x : x\in V(G) \}$.) We show
that $B_i = L_{\eta_i}^*V $ belongs to $\flgee$ by showing $R_fB_i
- B_i R_f = 0$ for all $f\in\edges(G)$ and $R_x B_i - B_i R_x = 0$
for all vertices $x \in V(G)$. From the previous lemma, if
$s(f)=x$, then
\[
R_f^* L_\eta = \big( R_f^* L_\eta \xi_x \big)\otimes \xi_x^* + L_\eta R_f^*,
\]
and so
\[
R_f L_\eta^* = - \xi_x \otimes \big( R_f^* L_\eta \xi_x\big)^* + L_\eta^* R_f.
\]
Consequently,
\begin{eqnarray*}
R_fB_i - B_i R_f &=& R_f L_{\eta_i}^* V - L_{\eta_i}^* V R_f
= \big( R_f L_{\eta_i}^* - L_{\eta_i}^* R_f \big) \,V \\
&=& \Big( - \xi_x \otimes \big( R_f^* L_{\eta_i} \xi_x \big)^* \Big) \,V
= - \xi_x \otimes \big( V^* R_f^* L_{\eta_i} \xi_x \big)^*.
\end{eqnarray*}

As each $\eta_i$ belongs to $ \M_x = Q_x \M$ for some
$x$, we have $ L_{\eta_i} \xi_x = 0$ if $\eta_i$ is not in $\M_x$. On the other hand, if $\eta_i
\in
\M_x$, then $L_{\eta_i} \xi_x = \eta_i$ and  this vector is orthogonal to $R_f V\H_G$
for each edge $f$. Hence $V^* R_f^* L_{\eta_i} \xi_x =0$ in this case. So we obtain
$
R_fB_i - B_i R_f = 0
$
for all $f$. Similarly, $R_x B_i - B_i R_x = 0$ for all $x$ since $B_i = L_{\eta_i}^* V$ is a
product of operators which commute with the projections $R_x = Q_x$. It follows that $B_i \in
\flgee$.

We now have $L_{\eta_i}^* V = B_i$ with $B_i \in \flgee$. Also,
since the range of $L_{\eta_i}$ is contained in the range of $V$,
it follows that $B_i$ is a partial isometry with range equal to
the range of $L_{\eta_i}^* $. But this range is the initial space
of $L_{\eta_i}$, which is easily seen to be $P_x$ when $\eta_i \in
\M_x$. Thus $P_x B_i = B_i\in\flgee$ is a partial isometry with
$\xi_x \in \ran (B_i)$. We claim that $B_i  = \lambda P_x$ for
some $|\lambda| = 1$. Indeed, there is a vector $\xi = B_i^* B_i
\xi$, $||\xi|| = 1$, such that $B_i \xi = \xi_x$ and by
considering the Fourier expansion for $B_i\in\flgee$ we can see
that $ \xi = B_i^* \xi_x = \ol{\lambda} \, \xi_x$ where $\lambda =
(B_i\xi_x, \xi_x)$.

Therefore $B_i \xi_x = \lambda \xi_x$ so that  $B_i \xi_w = R_w
B_i \xi_x = \lambda \xi_w$ for $w= xw$ and hence $B_i P_x =
\lambda P_x$. But  $B_i P_y = P_x B_i P_y = 0$ for $y\neq x$,
because $B_i$ is a partial isometry and otherwise we would have
$|| B_i^* P_x || > 1$  (a contradiction since $B_i^* P_x =
\ol{\lambda} P_x$). The claim now follows because
\[
B_i = B_i I = B_i\Big( \sum_{y\in V(G)} P_y\Big) = B_i P_x =
\lambda P_x.
\]

Evidently  $\big( L_{\eta_i} L_{\eta_i}^* \big) V = \lambda
L_{\eta_i} = L_{\eta_i^\prime}$ where $\eta_i^\prime = \lambda
\eta_i$. As the projections $L_{\eta_i} L_{\eta_i}^* $ are
orthogonal and sum to $VV^*$ it follows that
\[
V= VV^*V = \sum_i L_{\eta_i} L_{\eta_i}^* V = \sum_i L_{\eta_i^\prime}.
\]
Finally, each vector $\eta_i^\prime$ is supported on some $Q_x\H_G$ from the Beurling
Theorem.
But $L_{\eta_i^\prime} \xi_x = Q_x \eta_i^\prime$ by definition. Thus, as $V = \sum_i
L_{\eta_i^\prime}$ is a partial isometry it follows that the vectors $\{ \eta_i^\prime \}$ are
supported on distinct $Q_x \H_G$.
\bx

As an immediate consequence we obtain the following simple description of initial projections.
This result will be useful in Section~10.

\begin{cor}\label{stdformcor}
If $V$ is a partial isometry in $\flgee$, then the initial projection of $V$ is given by
\[
V^* V = \sum_{x\in \I} P_x
\qwhere \I = \big\{x\in V(G) :  V\xi_x \neq 0\big\}.
\]
\end{cor}

Lastly, we obtain an inner-outer factorization for elements of
$\flgee$ which generalizes the  $H^\infty$
\cite{Douglas_text,Hoffman_text} and  $\fL_n$ cases
\cite{DP1,Pop_beur}. Given a subset $\S \subseteq V(G)$ of
vertices,  define the {\it $\S$-inner} elements of $\flgee$ to be
the partial isometries with initial projection $\sum_{x\in\S}
P_x$. Also define the {\it $\S$-outer} elements to be those
elements of $\flgee$ with  range dense inside $\sum_{x\in\S}
P_x\H_G$.

\begin{cor}\label{innerouter}
Every $A$ in $\flgee$ factors as $A = VB$ where $V$ is an $\S$-inner element of $\flgee$ and
$B$ is $\S$-outer inside $\flgee$ with  $\S = \{ x\in
V(G) : A \xi_x \neq 0\}$.
\end{cor}

\Prf
Let $\M = \ol{\ran(A)} = \ol{A \frgee \xi_\phi} = \ol{\frgee A \xi_\phi}$ and let $\S = \{ x\in
V(G) : A \xi_x \neq 0\}$. Then there are unit $\frgee$-wandering vectors $\eta_x = Q_x
\eta_x$ for $x\in \S$ such that $\M = \sum_{x\in\S} \oplus L_{\eta_x} \H_G$. Let $V =
\sum_{x\in\S} \oplus L_{\eta_x}$ and observe that $V^* V = \sum_{x\in\S}
L_{\eta_x}^* L_{\eta_x} = \sum_{x\in\S} P_x$. Let $B_x = L_{\eta_x}^* A$ for $x\in\S$
and put $B =V^* A = \sum_{x\in\S} B_x$. It is clear that $A = VB$, and that each $B_x $ has
dense
range in $P_x \H_G$. Further, since the $L_{\eta_x}$ have pairwise orthogonal ranges which
span $\M$, it follows that $B$ has  range dense inside $\sum_{x\in\S} P_x \H_G$.
To complete the proof it suffices to show that each $B_x$ belongs to $\frgee^\prime = \flgee$,
and for this we may employ Lemma~\ref{stdformlemma} as in the proof of
Theorem~\ref{stdformthm}.

We mention there is also a uniqueness associated with the factorization $A =VB =
\sum_{x\in\S} L_{\eta_x} B_x$. The factors $L_{\eta_x}$ and $B_x$ of $A_x = L_{\eta_x}
B_x$ are unique up to a scalar
multiple since the wandering vectors $\eta\in Q_x \big( \M \ominus \sum_e \oplus R_e \M
\big)$, $x\in\S$, are unique up to a scalar.
\bx

%%%%%%%%%%%%%%%%%%%%%%%%%%%%%%%%%%%%%%%%%%%%%%
%%%%%%%%%%%%%%
\section{Classification and Automorphisms}\label{S:classification}
%%%%%%%%%%%%%%%%%%%%%%%%%%%%%%%%%%%%%%%%%%%%%%
%%%%%%%%%%%%%%

In this section we establish a classification theorem for the algebras
$\flgee$ by showing that
$G$ is a complete unitary invariant for $\flgee$. Our analysis also yields a large class of unitarily
implemented automorphisms of the algebras which act transitively on sets of eigenvectors of
$\flgee^*$.

%We say that two tree top graphs $G$ and $G^\prime$ are {\it equivalent} if their associated
%directed graphs are equivalent as directed graphs. This is the strongest notion of equivalence
%one
%could ask for, literally
%requiring the structure of the directed graphs to be identical, up to a relabelling of vertices and
%edges. For the sake of brevity we focus on the finite graph case, although we believe our
%theorem
%goes through for infinite graphs.

\begin{thm}\label{classification}
Let $G$ and $G^\prime$ be countable directed graphs. Then the following
assertions are equivalent.
\begin{itemize}
\item[$(i)$]$G$ and $G^\prime$ are isomorphic.
\item[$(ii)$] $\flgee$ and $\fL_{G^\prime}$ are unitarily equivalent.
\end{itemize}
\end{thm}

The proof relies on properties of $\flgee$ which are interesting in their own right.
We begin by showing that the family of vertex projections is a  unitary
invariant.
For $k \geq 0$, recall that $E_k$ is the projection onto $\spn \{
\xi_w : |w| =k\}$. Hence $E_0$ is the projection onto the vacuum
space $\spn \{ \xi_x : x\in V(G) \}$, and it is clear that
\begin{eqnarray}\label{cktidentity}
E_0 = \sum_{x\in V(G)} \xi_x\otimes \xi_x^* &=& I - \sum_{e\in \edges(G)}
L_e L_e^*.
\end{eqnarray}
Hence the rank one projections $\xi_x \otimes \xi_x^*$ are obtained by
compressing to $Q_x \H_G$
\begin{eqnarray}\label{vacuum}
\xi_x \otimes \xi_x^* = Q_x - \sum_{e\in\edges(G)}L_eL_e^* Q_x,
\end{eqnarray}
since the projections $\{ Q_x: x\in V(G) \}$ are reducing for
$\flgee$. The families $\{ R_e :
e\in\edges(G)\}$ and $\{ P_x : x\in V(G) \}$ have a similar
relationship.

\begin{lem}\label{minprojns}
The projections $\{ Q_x : x\in V(G) \}$ form the unique maximal
family of non-zero pairwise orthogonal irreducible projections for
$\fL_G$. The projections $\{ P_x: x\in V(G)\}$ play the same role for
$\frgee$.
\end{lem}

\Prf
As $I = \sum_x \oplus Q_x$, the projections $\{ Q_x \}$ form a maximal
family of pairwise orthogonal reducing projections for $\flgee$. To see minimality suppose
$0 \neq Q \leq Q_x$ is an $\fL_G$-reducing projection. From equation
(\ref{vacuum}) above, $Q$ commutes with $\xi_x\otimes \xi_x^*$ so the
vector $\xi_x$ either
belongs to $Q\H_G$
or is orthogonal to it. However, by the $\fL_G^*$-invariance of $Q \neq 0$
there is clearly a $\xi =
Q \xi$ such that $(\xi,\xi_x)\neq 0$.  Hence $\xi_x \in Q\H_G$, and
$\fL_G$-invariance gives
$Q=Q_x$. Thus each $Q_x$ is an irreducible projection.

To observe uniqueness of the family, suppose projections $\{Q_j\}$ form
another maximal such family.
The $Q_j$ are non-zero projections in $\flgee^\prime = \frgee$, hence by
Corollary~\ref{normals} each $Q_j$ belongs to the linear span of the
family $\{ Q_x \}$. Thus $Q_j$ is equal to the sum of a subset of these
projections, and
by the irreducibility of $Q_j$, it follows  that in fact $Q_j =
Q_x$ for some $x$.
Maximality of the $Q_j$ family ensures every $Q_x$ is obtained in this
manner. It
follows that the family of projections $Q_j$ must actually {\it be}
the family of $Q_x$. The same proof works for $\frgee$ and the family $\{
P_x \}$.
\bx

Our next step is to show how the number of directed edges between pairs of
vertices in $G$ can be computed in terms of $\flgee$ and the vacuum
vectors. Let $\flgee^0$ be the
$\wot$-closed two-sided ideal of
$\flgee$ generated by the  $ L_e$;
\[
\flgee^0 = \left< L_e : e\in E(G) \right>.
\]
Consideration of
Fourier expansions in $\flgee$ shows that every $A \in
\flgee^0$ satisfies $(A\xi_x, \xi_x ) =0$ for each $x$, and
in fact this condition characterizes $\flgee^0$.
This ideal helps to identify the directed graph structure of $G$ from the algebra $\flgee$.

\begin{lem}\label{formula}
Let
$a_{yx}$ be the number of directed edges $e =
y e x$ in  $G$ from vertex $x$ to vertex $y$. Then
\begin{eqnarray}
a_{yx} = \rank (P_yE_1Q_x ) = \dim \Big[ P_y \big( \flgee^0 \xi_x \ominus
(\flgee^0)^2 \xi_x
\big) \Big].
\end{eqnarray}
Furthermore, the ideal $\flgee^0$ may be computed as the set
\[
\flgee^0 = \big\{ A\in\flgee \,\, \, : \,\,\, A^* \xi_x = 0 \qfor x\in
V(G) \big\}.
\]
\end{lem}

\Prf
The family $\{ P_y, E_1, Q_x \}$ is mutually commuting, hence the
projection
$P_y E_1Q_x$ acts on basis vectors by $P_y E_1 Q_x \xi_w = E_1
\xi_{y w x}$. Thus it is the
range projection for the subspace
\[
P_y E_1 Q_x \H_G = \spn \{ \xi_e : e=y e x \}.
\]
This yields the identity $a_{yx} = \rank (P_y E_1 Q_x) $.
From the preceding discussion, every $A\in \flgee^0$ has an expansion of
the form $A \sim
\sum_{|w| \geq 1} a_w L_w$, and acts on basis vectors by $A\xi_v =
\sum_{|w| \geq 1}
a_w \xi_{wv}$. It follows that the subspace
\[
\flgee^0 \xi_x = \spn \{ \xi_w : s(w)=x \qand |w|\geq 1 \},
\]
whereas
\[
(\flgee^0)^2 \xi_x = \spn \{ \xi_w : s(w)=x \qand |w| \geq 2 \}.
\]
Therefore we have the subspace equalities
\begin{eqnarray*}
P_y \big( \flgee^0 \xi_x \ominus (\flgee^0)^2 \xi_x \big) &=&
\spn \{ \xi_w: w=y wx \qand |w| =1\} \\
&=& P_y E_1 Q_x \H_G.
\end{eqnarray*}

To see the alternate description of $\flgee^0$ in terms of vacuum vectors,
notice that if $A \sim
\sum_{|w| \geq 1} a_w L_w$ belongs to $\flgee^0$, then clearly $A^* \xi_x
= 0 $ for
 $x\in V(G)$. On the other hand, if $A \sim \sum_w a_w L_w$ belongs to
$\flgee$ and $A^*$ annihilates each vacuum vector, then
\[
0 = (\xi_x , A^* \xi_x) = (A \xi_x , \xi_x) = a_{x} \qfor x\in V(G),
\]
and $A$ belongs to $\flgee^0$.
\bx

We now  prove the classification theorem. The proof makes use of Theorem~\ref{autothm},
which we present below because of its independent interest. But notice that this theorem is not
needed for the special case of graphs with no loop edges at vertices; in this case
Proposition~\ref{noevalues} shows that the ideal $\flgee^0$ is invariant under unitary
isomorphism.

\vspace{0.1in}

{\noindent}{\it Proof of Theorem~\ref{classification}.}
If  $G$ and $G^\prime$ are isomorphic, then there is a relabelling map
between the vertices and edges of the graphs which clearly induces a unitary operator between
$\H_G$ and $\H_{G^\prime}$. This unitary intertwines the generators $\{
L_e, P_x \}$ of $\flgee$ with the generators $\{L_{e^\prime}, P_x^\prime\}$ of
$\fL_{G^\prime}$, and {\it a fortiori} the algebras $\flgee$ and $\fL_{G^\prime}$ are unitarily
equivalent.

Conversely, suppose there is a unitary operator $U: \H_{G^\prime} \rightarrow
\H_G$ for which $U^* \flgee U = \fL_{G^\prime}$. Without loss of generality, we may assume
that $\H_G = \H_{G^\prime} $. By Lemma~\ref{minprojns} the family $\{P_x\}$ is a unitary
invariant of $\flgee$, hence the number of vertices in $G$ and
$G^\prime$ is the same and we
may assume, perhaps after reordering, that $U^* P_x U = P_x^\prime$ for each $x$ in $V(G) =
V(G^\prime)$. Further,
under this
unitary equivalence the vacuum vectors $\xi_x^\prime$ for $\fL_{G^\prime}$ are mapped to
eigenvectors
for $\fL_{G}^*$. But Theorem~\ref{autothm} shows that these vectors may be moved to the
vacuum vectors $\xi_x$ by a unitary which implements an automorphism of
$\fL_{G}$. In particular, we may also assume that $U \xi_x^\prime = \xi_x$ for each $x$.
It remains to show that the number $a_{yx}$ of directed edges $e=yex$ in $G$ is equal to the
number $a^\prime_{yx}$ of directed edges $e^\prime = ye^\prime x$ in $G^\prime$, for all
$x,y\in V(G) = V(G^\prime)$. This readily follows from Lemma~\ref{formula} and our above
assumptions; indeed, we have $U^* \flgee^0 U = \fL_{G^\prime}^0$ and
\begin{eqnarray*}
a_{yx}^\prime &=&  \dim \Big[ P_y^\prime \big( \fL_{G^\prime}^0 \xi_x^\prime \ominus
(\fL_{G^\prime}^0)^2 \xi_x^\prime
\big) \Big] \\
&=& \dim \Big[ U^* P_y \big( \flgee^0 \xi_x \ominus
(\flgee^0)^2 \xi_x
\big) \Big] = a_{yx}.
\end{eqnarray*}
Therefore the directed graphs $G$ and $G^\prime$ are isomorphic.
\bx

The proof of the classification theorem relies on the existence of
certain unitarily implemented automorphisms of $\flgee$ which we
now discuss. Let $\nu = \nu_{\lambda, x} = P_x Q_x \nu_{\lambda,
x}$ be an eigenvector for the algebra $\frgee^*$, as given in the
$\frgee$ version of Theorem~\ref{eigenvalues} (see
Remark~\ref{evectorremark}). We show there is an automorphism of
$\flgee$ which is implemented by a unitary which maps
$\nu_{\lambda, x}$ to $\nu_{0,x} = \xi_x$. This transitive action
of unitary automorphisms was obtained for free semigroup algebras
by Davidson and Pitts \cite{DP2} and our proof is an elaboration
and generalization of their analysis.

Recall first that in the case of the free semigroup algebra $\fL_n$ and its commutant algebra
$\fR_n$, with eigenvector $\eta$ for $\fR_n^*$, the $\fR_n$-wandering subspace $\M$ for the
$\fR_n$-invariant subspace $\{ \eta \}^\perp$ has an orthonormal basis consisting of $n$ vectors
$\{ \eta_1, \ldots, \eta_n \}$ say. Also recall that if a vector $\eta$ is $\fR_n$-wandering then one
can define the isometry $L_\eta$ in $\flgee$ by the specification $L_\eta \xi_w = R_w \eta$ for
all words $w$. The desired automorphism in this case is in fact effected by the correspondence
$L_{e_i} \rightarrow L_{\eta_i}$, for $1 \leq i \leq n$, where $e_1, \ldots, e_n$ are the $n$ loop
edges of the free semigroup graph. Using only the fact that the wandering
space for $\M$
is $n$-dimensional, we shall develop a similar argument from first principles to define the
desired automorphism of $\flgee$. There are some complications in the new setting. In
particular, to capture all the generators we must consider more than the wandering subspace for
$\{ \nu \}^\perp$.

For the next theorem we use the following notation. The set $W_x$ is the
collection of words in
the loop edges at vertex $x$, and $\H_x$ is the closed span of the basis vectors $\xi_w$ for $w$
in $W_x$. Also we identify $\H_x$ with the Fock space for $\fL_n$ and identify the generators
$L_1, \ldots, L_n$ (respectively $R_1, \ldots, R_n$) of $\fL_n$ (respectively $\fR_n$) with the
restrictions $L_e|_{\H_x}$ (respectively $R_e|_{\H_x}$) for the $n$ edges $e$ with $x e x =
e$. In particular, from our earlier notation we have $n = a_{xx}$.

\begin{thm}\label{autothm}
Let $G$ be a countable directed graph, and let $\nu = \nu_{\lambda, x}
=Q_x P_x \nu_{\lambda, x}$ be an eigenvector for
$\frgee^*$.
\begin{itemize}
\item[$(i)$] The subspace
\[
\M = \{\nu \}^\perp \bigcap \big( \cap_{y \neq x} \{ \xi_y \}^\perp \big)
\]
is $\frgee$-invariant with wandering subspace basis $\{ \eta_e : e
\in \edges(G) \}$ where $\eta_e =\xi_e$ if $r(e)\neq x$, where
$\eta_e = R_e \nu$ if $r(e)=x$ and $s(e)\neq x$, and where $\{
\eta_e : x e x = e \}$ is a basis for the $\fR_n$-wandering
subspace of
\[
\{\nu \}^\perp \cap \H_x
\qfor \{ R_j \} = \{ R_e|_{\H_x} : x e x = e \}.
\]
\item[$(ii)$] The correspondence $e \rightarrow L_{\eta_e}$, and $x
\rightarrow L_x$ gives a purely atomic free
partial isometry representation of $G$ satisfying the multiplicity one condition at each vertex,
and the correspondence $L_e \rightarrow L_{\eta_e}$ extends to an automorphism of
$\flgee$.
\item[$(iii)$] Let $S_e = L_{\eta_e}$, for $e\in \edges(G)$, and $S_y = L_y$ for $y \in V(G)$.
Then there is a unitary operator $W$ defined by
\[
W \xi_w = \left\{ \begin{array}{cl}
w(S) \nu & \mbox{if $s(w)=x$} \\
w(S) \xi_y & \mbox{if $s(w)=y$, $y\neq x$}.
\end{array}\right.
\]
Furthermore, $\Ad W (L_e) = W L_e W^* = S_e$ for $e \in \edges(G)$ and $\Ad W$ is the
automorphism of
$\flgee$ given in $(ii)$.
\end{itemize}
\end{thm}

\Prf
Choose $\eta_{e_1}, \ldots, \eta_{e_m}$  to
be an orthonormal basis (possibly countably infinite) for
the $\fR_n$-wandering subspace of $\{ \nu \}^\perp \cap \H_x$, and let $\eta_e$, for the other
edges of $G$, be specified as in $(i)$. We claim this basis spans the $\frgee$-wandering
subspace for $\M$ which is $\M \ominus \big( \sum_e \oplus R_e
\M \big)$.
Indeed, from the definition of $\M$ it is easy to see that each of the three types
of basis vector belong to this wandering subspace. On the other hand,
by the choice of $\eta_e$ for $e$ with $e=xex$ it is not hard to see that
all the
non-zero vectors of the form $R_w \eta_e$ give rise to an orthonormal
basis for $\M$. Thus the
wandering subspace for $\M$ has orthonormal basis $\{ \eta_e :
e\in\edges(G)\}$. Further, each $\eta_e$ is supported on a particular
$Q_x$, hence the decomposition $\M = \sum_{e\in\edges(G)} \oplus \frgee
[\eta_e]$ really is the decomposition of $\M$ into minimal cyclic
subspaces indicated in the $\frgee$ version of Theorem~\ref{beurlingthm} (see
Remark~\ref{beurlingremark}).

We establish $(ii)$ and $(iii)$ together. The operators $S_e =
L_{\eta_e}$ have the defining property $S_e \xi_w = R_w \eta_e$,
for all  $w$ in $\fngee$. It follows that $S_e$ is a partial
isometry with initial projection $P_y = L_y$ where $s(e)=y$ since
$\eta_e = Q_y \eta_e$. Moreover, the final projection of $S_e$ is
the space $\frgee [\eta_e]$ by design and so, by our construction
of the basis it is clear that
\begin{eqnarray*}
\sum_{r(e)=y} S_e S_e^* &=& P_y - \xi_y \otimes \xi_y^*, \qfor y\neq x, \\
\sum_{r(e)=x} S_e S_e^* &=& P_x - \nu \otimes \nu^*.
\end{eqnarray*}
It now follows that the map $e \rightarrow S_e$, $y \rightarrow S_y=L_y$ gives a
free partial
isometry representation of $G$ satisfying the multiplicity one condition
considered in Proposition~\ref{repntheory} and that the
generators $\{S_e, L_y \}$ are mutually unitarily equivalent to the generators $\{L_e, L_y \}$ by
the unitary $W$.

It remains to show that the unitary automorphism $\Ad W(X) = WXW^*$ satisfies $\Ad
W(\flgee) = \flgee$, that is, that the unitary automorphism $\Ad W$ of $\bofh$ restricts to an
automorphism of $\flgee$ rather than an endomorphism.
Using the  gauge automorphisms of
Section~\ref{S:freeprods}, at this point in the proof we can easily follow the
free semigroup approach  (see \cite{DP2}, Remark 4.13). Indeed,
the algebra $\fA = \Ad W(\flgee)$ is contained in $\flgee$, and hence
$\nu_{0,x} = \xi_x$ is an eigenvector for $\fA^*$. Since $\fA$ is
unitarily equivalent to $\flgee$, there is a non-zero $\mu$ such that $W
\nu_{\mu,x} = \xi_x$. Hence we can apply the above argument again to
obtain another unitary $W^\prime$ for which $\Ad W^\prime W (L_e) =
S_e^\prime$, where the $S_e^\prime$ are determined as above by an
orthonormal basis $\B$ for the wandering subspace of the subspace $\cap_x
\{\xi_x\}^\perp$. Let $U$ be a unitary in $\U_m$ which intertwines the
orthonormal set $\{ \xi_e : e \in \edges(G) \}$ with the vectors of $\B$,
in such a way that $UP_y = P_y$ for each $y$. Then it follows that
$\Ad W^\prime W = \theta_U$, the gauge automorphism of $\flgee$ determined
by $U$. Consequently, the two endomorphisms of $\flgee$ must actually be
automorphisms.
\bx

We can immediately deduce the corresponding classification of the norm-closed
algebras generated by the generators of $\flgee$. Let us denote this algebra,
which is a non-commutative version of the disc algebra, as $\A_G$. In the case of
a finite directed graph this algebra was studied in the general framework of
tensor algebras over correspondences by Muhly and Solel
\cite{Muhly_fd,MS_CJM}, but the basic classification question was not
considered.
Recall that $\A_G$ and $\A_{G^\prime}$ are star-extendibly isomorphic if there is an
isomorphism $\A_G \rightarrow \A_{G^\prime}$ which is the restriction of a (necessarily
unique) $\ca$-algebra isomorphism $\ca(\A_G) \rightarrow \ca(\A_{G^\prime})$.

\begin{cor}\label{normclass}
Let $G$, $G^\prime$ be countable directed graphs. Then the following assertions are
equivalent.
\begin{itemize}
\item[$(i)$] $G$ and $G^\prime$ are isomorphic.
\item[$(ii)$] $\A_G$ and $\A_{G^\prime}$ are unitarily equivalent.
\end{itemize}
Also, if each vertex of $G$, $G^\prime$ has finite degree then a star-extendible
isomorphism between $\A_G$ and $\A_H$ is unitarily implemented.
\end{cor}

\Prf
If $\A_G$ and $\A_{G^\prime}$ are unitarily equivalent, then so are their $\wot$-closures
and so the equivalence of $(i)$ and $(ii)$ follows from
Theorem~\ref{classification}. For the final assertion note that  $\ca(\A_G)$
contains the collection of compact operators
$K\in\bofh$ such that $KQ_x = Q_x K$ for all vertices $x$.
Indeed, for $x \in V(G)$ and words $v=vx$, $w=wx$ in $\bbF^+(G)$, equation
(\ref{cktidentity}) shows that the rank one projection
\begin{eqnarray*}
\xi_v \otimes \xi_w^* &=& L_v \left( \xi_x \otimes \xi_x^* \right) L_w \\
&=& L_v P_x \left( I - \sum_{e\in\edges(G)} L_e L_e^* \right) P_x L_w\\
&=& L_v P_x L_w - \sum_{e\in\edges(G)} L_v P_x L_e L_e^* P_x L_w
\end{eqnarray*}
belongs to $\ca (\A_G)$.
Note that the summation here is finite.
Thus it follows that the isomorphism
$\ca(\A_G)\rightarrow\ca(\A_{G^\prime})$ is unitarily implemented, and hence the restriction
to $\A_G$
produces a unitary equivalence with $\A_{G^\prime}$.
\bx

We finish this section with an example which may clarify some of the subtleties of
Theorem~\ref{autothm}.

\begin{eg}\label{autoeg}
 Let $G$ be the  directed graph with vertex set $V(G) = \{ x_1, x_2\}$ and
edges  $e_i = x_1 e_i x_1$ for $i = 1,2$, $e_3 = x_2e_3x_1$, and
$e_4 = x_1 e_4 x_2$. Let $\nu = \nu_{\lambda, 1}= P_1 Q_1 \nu_{\lambda,
1}$ be an eigenvector for $\frgee^*$. In this case
$
\M = \{ \nu \}^\perp \bigcap \{ \xi_2 \}^\perp.
$
The orthonormal basis for the wandering space of $\M$ from the theorem is
given by
$\eta_{e_3} = \xi_{e_3}$ since $e_3 = x_2e_3$, $\eta_{e_4} = R_{e_4} \nu$ since $e_4 = x_1
e_4 x_2$, and $\{ \eta_{e_1}, \eta_{e_2} \}$ is a basis for the $\fR_2$-wandering subspace of
$\{ \nu \}^\perp \cap \H_1$, where $\H_1 = \spn \{ \xi_w : w\in W_1 \}$
is identified with the Fock space for $\fL_2$ and $\fR_2$, and
$W_1$ the set of words in $e_1$, $e_2$.
Thus, as in the theorem we have $\M = \sum_{i=1}^5 \oplus \frgee
[\eta_{e_i}]$.

The basis $\B$ in the proof will form an orthonormal basis for the
wandering space of $\{\xi_1\}^\perp \cap \{ \xi_2\}^\perp$, which
is $\spn \{ \xi_e : e\in\edges(G)\}$. From the construction
outlined in the statement of the theorem, this basis will also
have each of its vectors fully supported on some $P_j$. A gauge
unitary $\td{U}$ of the type used in the proof will be determined
here by a unitary $U \in \U_4$ which fixes $\xi_{e_3}$ and
$\xi_{e_4}$ and is allowed to scramble the subspace
$\spn\{\xi_{e_1}, \xi_{e_2}\}$.
\end{eg}

%%%%%%%%%%%%%%%%%%%%%%%%%%%%%%%%%%%%%%%%
\section{Partly Free Algebras}\label{S:partfree}
%%%%%%%%%%%%%%%%%%%%%%%%%%%%%%%%%%%%%%%%

We now determine in graph-theoretic terms when an operator algebra $\flgee$
contains the free semigroup algebra $\fL_2$ as a subalgebra. More generally, let us say that a
$\wot$-closed operator algebra
$\fA$ is {\it partly free} if it contains the free semigroup algebra $\fL_2$ as a subalgebra in the
sense of the following definition.
\begin{defn}\label{partlyfreedefn}
A $\wot$-closed algebra $\fA$ is {\it partly free} if there is an inclusion map $\fL_2
\hookrightarrow \fA$ which is the restriction of an injection between the generated von
Neumann algebras. If the map can be chosen to be unital, then $\fA$ is said to be {\it unitally
partly free}.
\end{defn}
These notions parallel somewhat the requirement that a $\ca$-algebra contain
$\O_2$, or that a discrete group contain a free group.
Theorems~\ref{doublecyclethm} and \ref{strongdoublecyclethm} determine when the
algebras $\flgee$ are partly free and unitally partly free. We first set aside two results which
have  intrinsic interest.

\begin{lem}\label{doublelemma1}
Let $\eta \in Q_x\H_G$ be a unit $\frgee$-wandering vector.
Suppose that $L_\eta L_\eta^*\leq
L_\eta^* L_\eta$. Then $\eta = \sum_u a_u \xi_u$ and for each $u\neq xux$,
$a_u = 0$. That is, $\eta$ is supported on basis vectors corresponding to
words forming cycles at the vertex $x$ in $G$.
\end{lem}

\Prf
Since $\eta$ belongs to $Q_x \H_G$, it follows that $L_\eta^* L_\eta = P_x$, whence $\eta =
L_\eta \xi_x \in P_x$ by assumption.
Further, the non-zero vectors among $ L_\eta \xi_u $ form an orthonormal set, hence
\[
1 = \sum_u |a_u|^2 = ||\eta||^2
= ||L_\eta \eta ||^2 = \sum_u |a_u|^2 ||L_\eta \xi_u||^2.
\]
Thus if $a_u\neq 0$ then $||L_\eta\xi_u||^2 = 1$. In particular,
$\xi_u$ belongs to the initial space of $L_\eta$ which is
$L_\eta^* L_\eta\H_G = P_x \H_G$. Thus $r(u)=x$, but since $\eta =
\sum_u a_u \xi_u$ is in $Q_x\H_G$ we also have $s(u)=x$ when $a_u
\neq 0$. It follows that $u=xux$ for all $u$ with $a_u \neq 0$.
\bx

The cycle algebras $\fL_{C_n}$ from Example~\ref{cyclegraphs} give the
motivational
subclass of infinite-dimensional  algebras which are not partly free.

\begin{lem}\label{doublelemma2}
The cycle algebras $\fL_{C_n}$, $1 \leq n < \infty$, do not contain pairs of partial
isometries $U$, $V$ which satisfy condition $(iii)$ of
Theorem~\ref{doublecyclethm}.
\end{lem}

\Prf
This readily follows from the  matrix function theory description of the cycle
algebras $\fL_{C_n}$ since
a similar fact holds in the algebras $H^\infty \otimes \M_n$ and their
direct sums. This in turn follows from elementary Toeplitz operator theory, or from the fact that
these algebras possess a natural faithful trace.
\bx

We now define the graph-theoretic notions we require. Recall that a cycle in a directed graph is
{\it minimal} if it is not a power of another cycle.

\begin{defn}\label{doubledefn}
We say $G$ has the {\it double-cycle property} if there
are distinct minimal cycles $w=xwx$, $w^\prime =xw^\prime x$ over the same vertex $x$ in
$G$. We say $G$ has the {\it strong double-cycle property} if for every vertex $x$ in $G$ there
is a directed path from $x$ to a vertex lying on a double-cycle.
\end{defn}

\begin{thm}\label{doublecyclethm}
The following assertions are equivalent for a countable directed graph $G$ with a finite number
of vertices.
\begin{itemize}
\item[$(i)$] $G$ has the double-cycle property.
\item[$(ii)$] $\flgee$ is partly free.
\item[$(iii)$] There are non-zero partial isometries $U$, $V$ in $\flgee$ with
\[
U^*U = V^*V, \,\,\,\, UU^*\leq U^*U, \,\,\,\, VV^*\leq V^*V, \,\,\,\, U^* V =0.
\]
\end{itemize}
\end{thm}

\Prf
For $(i)\Rightarrow (ii)$, observe that   if $w$,
$w^\prime$ are cycles of minimal length at vertex $x$, and $w\neq
w^\prime$, then we
may take $U= L_w$, $V= L_{w^\prime}$ to define an injection of $\fL_2$ into $\flgee$. Since
$(ii)$ clearly implies
$(iii)$, it remains to establish the implication $(iii) \Rightarrow (i)$.

By Theorem~\ref{stdformthm} we have the initial projection $U^*U = V^*V$
equal to the sum of certain $P_i\equiv P_{x_i}$. Without loss of generality let us
assume
\begin{eqnarray}\label{induction}
U^*U = V^*V = P_1 + \ldots + P_k.
\end{eqnarray}
We establish $(i)$ by induction on $k$.

For $k=1$, observe that Lemma~\ref{doublelemma1} gives a double-cycle
over $x_1$ when $U^*U
= V^*V = P_1$. Indeed, suppose by way of contradiction, that $G$ fails to
have the double-cycle property. As $k=1$, Theorem~\ref{stdformthm} gives $U=
L_\eta$ and $V=L_{\eta^\prime}$.  (Observe that there is at least a single loop
edge over $x_1$ since $U$ and $V$ are non-zero.)
By Lemma~\ref{doublelemma1} we deduce that for some minimal cycle $w$ (possibly a
single loop edge) both $\eta$ and $\eta^\prime$ belong to the subspace
\[
\H_w = \spn \{ \xi_x, \xi_{w^m} : m=1,2, \ldots  \}.
\]
But  $\H_w$ can be identified with $H^2$, and $L_w$ is then identified
with the unilateral shift on $H^2$. Consider the subspaces
\[
H_\eta = \spn \{ R_w^m \eta : m\geq 0 \} \qand
H_{\eta^\prime} = \spn \{ R_w^m \eta^\prime : m\geq 0 \}.
\]
Since these are non-zero invariant subspaces for the multiplicity-one
unilateral shift $L_w|_{\H_w}$, it follows from the classical Beurling
theorem for $H^2$ that they have non-empty intersection. This contradicts
the hypothesis, since $\H_\eta \subseteq \ran (L_\eta)$ and
$\H_{\eta^\prime} \subseteq \ran (L_{\eta^\prime})$.

Let $k\geq 2$ and assume  $(iii)\Rightarrow (i)$ holds for
$m=1,\ldots
,k-1$; that is, $G$ contains a double-cycle whenever $\flgee$ contains a
$U$,$V$ satisfying $(iii)$ for which $U^*U=V^*V$ is a sum of at most $k-1$
projections $P_i$. Let  $\S = \{ x_1, \ldots, x_k\}$ be the vertices corresponding to the
projections $P_i$ in (\ref{induction}).  We may assume that every vertex $x\in \S$
has the property that a
directed path in $G$  leaves it for another vertex in
$\S$. For if $x\in \S$ was a vertex without this property, then  $UP_x = P_x UP_x$ and $VP_x =
P_x VP_x$ would be non-zero  partial isometries
in $\flgee$ with pairwise orthogonal ranges and  initial projection $P_x$
containing their final projections. Thus, by the $k=1$ case, $G$ would contain a double-cycle.

Now fix $x\in \S$ for the moment and consider a directed path $w$ in $G$ that has initial vertex
$x$ and final vertex in $\S$, and
passes through a
maximal number of vertices in $\S$ without going through the same vertex in $\S$ twice. Let $y
\in \S$ be the final
vertex of $w = ywx$. We know there is a path from $y$ to another vertex
$z$ in $\S$, but by maximality $w$ must pass through $z$. Consequently,
there is a subset $\A \subseteq \S$ of vertices which lie on a cycle.
Let us  assume this cycle does not
cross itself, and further assume there are no paths in $G$ outside the cycle which have both
initial and
final vertices belonging to the set of vertices which form the cycle
(otherwise $G$ would clearly contain a double-cycle).

Then $\A$ is a proper subset of $\S$. To see
this, suppose $\A = \S$, and let $P$ be the projection which is the sum of
all $P_x$ for which $x$ is a vertex on the cycle. Then by the assumptions on $\A$ in the
previous
paragraph, the
algebra $P \flgee P$ will consist of operators in $\flgee$ which have
non-zero Fourier coefficients only for basis vectors
corresponding to words whose letters are edges in the cycle.
Let $P_0 = \sum_{x\in \A} P_x \leq P$. Then evidently $P \flgee|_{P\H_G}$ is
unitarily
equivalent to $\fL_{C_n}$ for some $n$, and $P U|_{P\H_G}= P_0 UP_0|_{P\H_G}=
U|_{P\H_G}$, $P V|_{P\H_G}= P_0 VP_0|_{P\H_G}= V|_{P\H_G}$ would yield a pair of
partial isometries in $\fL_{C_n}$ satisfying
condition $(iii)$. But this cannot happen by
Lemma~\ref{doublelemma2}. Thus $\A$ must in
fact be a proper subset of $\S$.

Let $\B$, $\C$, $\D$ be the subsets of $\S$ which make up the complement
of $\A$ consisting of respectively:
final vertices of paths with initial vertices in $\A$; vertices for which there is
a path that leaves it and ends at a vertex in $\A$; and vertices in $\S$ for
which there are no paths to or from vertices in $\A$.
Thus $\S\setminus \A = \B\cup\C\cup\D$.
We can assume that $\B$ is non-empty.
For otherwise, there would be no  edges which emerge from the cycle graph of $\A$ and the
above reduction argument can be applied, together with the fact that the cycle
algebras are not partly free, to view $U$ and $V$ as elements of $\fL_H$ where $H$
is the graph obtained when the sink vertex set $\B$ is removed. But from the
definition of $\D$, there are no directed paths from a vertex in
$\B$ to a vertex in $\D$. Further, there are no paths from $\B$ to $\C$ by
the assumptions on $\A$. Thus there are no paths in $G$ from $\B$ to {\it
any} of $\A$, $\C$, or $\D$. Let $P = \sum_{x\in\B} P_x$ be the sum of
projections corresponding to vertices in $\B$. Then $UP = PUP$, $VP =PVP$
are non-zero and $P$ is the sum of strictly less than $k$ of the $P_x$'s.
Hence by induction $G$ has a double-cycle.

Therefore we conclude that $G$ does indeed contain a double-cycle when
condition $(iii)$ holds, and this completes the proof.
\bx

We next establish the unital version of the previous theorem.

\begin{thm}\label{strongdoublecyclethm}
The following assertions are equivalent for a countable directed graph $G$ with a finite number
of vertices.
\begin{itemize}
\item[$(i)$] $G$ has the strong double-cycle property.
\item[$(ii)$] $\flgee$ is unitally partly free.
\item[$(iii)$] There are  isometries $U$, $V$ in $\flgee$ with
\[
U^* V =0.
\]
\end{itemize}
\end{thm}

\Prf
Condition $(ii)$ clearly implies $(iii)$. For $(iii)\Rightarrow (i)$,
notice that in the proof of Theorem~\ref{doublecyclethm} we actually showed that from every
vertex $x\in V(G)$ with $P_x \leq U^* U = V^*V$, there is a directed path into a double-cycle.
Thus, in this case, we may apply this argument to every vertex in $G$ since
$U^*U =V^*V = I = \sum_{x\in V(G)} P_x$. In particular, $G$ satisfies the strong double-cycle
property when $(iii)$ holds.

We next establish  $(i)\Rightarrow(iii)$ and  $(i)\Rightarrow (ii)$ together.  Thus suppose $G$
satisfies the strong double-cycle property.
Fix a double-cycle
in $G$ and let $\B$ be the (maximal) collection of all vertices which lie on paths going into or
on this double-cycle. Let $x$ be a vertex in $\B$ which belongs to the given double-cycle. Then
there are minimal cycles $w_1 = xw_1x \neq w_2 = xw_2x$. Let
$\bbF_2^+(w_1,w_2)$ be the set of all words in the generators $w_1$, $w_2$ and consider the
subspace
\[
\H_{w_1,w_2} = \spn \big\{ \xi_x, \xi_w: w\in \bbF_2^+(w_1,w_2) \big\}.
\]
Fix a positive integer $k$ such that $2^k \geq 2 |\B|$. Amongst the $2^k$ words of length $k$ in
$\bbF_2^+(w_1,w_2)$, choose a set of cardinality $2|\B|$ and label elements of this set by
$\{u_{y}^{(i)}:y\in\B, i=1,2\}$. For every $y\in\B$ there is a path $v_y$ such that $v_y =
xv_yy$. For $y\in\B$ and $i=1,2$ let $w_y^{(i)}$ be the path
$w_{y}^{(i)} = u_{y}^{(i)} v_y$.
Observe that each of the partial isometries $L_{w_y^{(i)}}$ has initial projection
$L_{w_y^{(i)}}^* L_{w_y^{(i)}} = P_y$. Further, the entire family of operators
$\{ L_{w_y^{(i)}}: y\in\B, i=1,2\}$ have pairwise orthogonal ranges by design. Thus it follows
that the operators
\[
U_\B = \sum_{y\in\B}\oplus L_{w_{y}^{(1)}} \qand V_\B = \sum_{y\in\B} \oplus
L_{w_{y}^{(2)}}
\]
are partial isometries in $\flgee$ with mutually orthogonal ranges contained in $\H_{w_1,w_2}$
and initial projections satisfying
\[
U^*_\B \, U_\B = \sum_{y\in\B} P_y  = V^*_\B \,V_\B.
\]

Now let $\B_1, \ldots, \B_d$ be a maximal family of disjoint sets of vertices of
$G$, where each of
these sets is obtained in the same manner as the above set $\B$. (Choose $\B_1$ as $\B$ was
chosen, then obtain $\B_2$ in a similar manner from $V(G)\setminus \B_1$, et cetera.) Since
the strong double-cycle property
holds for $G$, the disjoint union $\cup_i \B_i = V(G)$.
Let $U_{\B_1}, \ldots, U_{\B_d}$ and $V_{\B_1}, \ldots, V_{\B_d}$ be partial isometries
obtained as in the construction of the previous paragraph.
Then the operators $\{ U_{\B_i}, V_{\B_j} : 1\leq i,j \leq d\}$ have
pairwise orthogonal ranges with initial projections satisfying
\[
U_{\B_i}^* U_{\B_i} = \sum_{y\in\B_i} P_y = V_{\B_i}^* V_{\B_i} \qfor 1\leq i \leq d.
\]
Therefore it follows that the operators
$U = \sum_{i=1}^d \oplus U_{\B_i}$ and $V = \sum_{i=1}^d \oplus V_{\B_i}$
are isometries in $\flgee$ which have mutually orthogonal ranges, and hence condition
$(iii)$ holds.
Finally, the map which sends the two generators of $\fL_2$ to $U$ and $V$ induces an injection
of $\fL_2$ into $\flgee$, and $(ii)$ holds. This completes the proof.
\bx

\begin{rem}
In the finite graph case
it is clear from the proof  of Theorem~\ref{strongdoublecyclethm} that
the family of paths which determine the partial isometries $L_w$ in the sums defining $U$ and
$V$ can be chosen so that they all have the same length. Hence it
follows that $G$ has the strong double-cycle property precisely when the transpose graph $G^t$
satisfies the {\it entrance condition} from \cite{MS_CJM} (c.f. Definition~5.8), used  as a
condition
which guarantees the existence of isometries with mutually orthogonal ranges in the commutant.
Thus in the finite graph case of
Theorem~\ref{strongdoublecyclethm} we have proved this entrance condition on $G^t$ is
actually {\it equivalent} to the existence of isometries with mutually orthogonal ranges in the
commutant.
\end{rem}

We finish with a brief discussion of hyper-reflexivity. Given an operator algebra
$\fA$, a measure of the
distance to $\fA$ is given by
\[
\beta_{\fA} (X) = \sup_{L\in\Lat\fA} ||P_L^\perp X P_L||,
\]
where $P_L$ is the projection onto the subspace $L$ and $\Lat \fA$ is the lattice of invariant
subspaces for $\fA$. Evidently,
$\beta_\fA (X) \leq \dist (X,\fA)$, and the algebra $\fA$ is said to be
{\it hyper-reflexive} if there is a constant $C$ such that $\dist(X,\fA)
\leq C \beta_\fA (X) $ for all $X$.

The list of known hyper-reflexive algebras is short, but growing. See
\cite{Arv_nest,Berc,Chr,Dav_Toep,DP1} for examples appearing in the
literature. For the algebras $\fL_n$, hyper-reflexivity was
established by Davidson for $\fL_1 = H^\infty$ \cite{Dav_Toep}, and by
Davidson and Pitts for the free semigroup algebras $n\geq 2$ \cite{DP1}. In \cite{Berc}
Bercovici introduced a
general method motivated by the $\fL_n$ case, and lowered the upper bound for the $\fL_n$
distant constant. In particular, he proved that an algebra is hyper-reflexive
with distant constant no greater than $3$ whenever its commutant contains
a pair of isometries with orthogonal ranges.

\begin{cor}\label{hyperreflex}
Let $G$ be a countable directed graph with finitely many vertices for which the transpose graph
$G^t$
satisfies the strong double-cycle property.
Then $\flgee$ is hyper-reflexive with distant constant at most 3.
\end{cor}

\Prf
{}From Lemma~\ref{commutantlemma} the commutant $\flgee^\prime = \frgee$ is unitarily
equivalent to
$\fL_{G^t}$, which contains a pair of isometries with pairwise orthogonal ranges by the
previous theorem.
Thus the result follows from a direct application of Bercovici's result.
\bx

\begin{rem}
Using Corollary~\ref{hyperreflex} and separate arguments for graphs without the double cycle
property in the transpose graph it can be shown that $\flgee$ is hyper-reflexive for every finite
graph \cite{JP}. It would be interesting to have a characterization of general `hyper-reflexive
graphs', although this is likely to be a deep problem.
\end{rem}

%%%%%%%%%%%%%%%%%%%%%%%%%%%%%%%%%%%%
%\section{Memo Section}
%%%%%%%%%%%%%%%%%%%%%%%%%%%%%%%%%%%%

%%%%%%%%%%%%%%%%%%%%%%%%%%%%%%%%%%%%%%%%%%%%%%
%%%%%%%%%%%%%
%\section{$\wot$-closed Ideals of $\flgee$}\label{S:ideals}
%%%%%%%%%%%%%%%%%%%%%%%%%%%%%%%%%%%%%%%%%%%%%%
%%%%%%%%%%%%%

%\begin{note}
%In this section: Lattice iso thm between right and two-sided ideals, and
%the corresponding subspace lattices?? note on what goes wrong for the left
%ideals. Characterization of all mult linear fnls, from
%eigenvectors?? Characterization of commutator ideal, obvious
%thing?? Discuss factorization in right ideals, will have to adjust in
%infinite graph case.
%\end{note}

%%%%%%%%%%%%%%%%%%%%%%%%%%%%%%%%%%%%%%%%%%%%%%
%%%%%%%%%%%%%
%%%%%%%%%%%%%%%%%%%%%%%%%%%%%%%%%%%%%%%%%%%%%%
%%%%%%%%%%%%%

{\noindent}{\it Acknowledgements.} We would like to thank the
referees for a number of helpful suggestions. We are grateful to
Paul Muhly for pointing out connections with his work and to
Frederic Jaeck for discussions on matrix function algebras. We
also thank Elias Katsoulis for detecting a gap in the proof of
Theorem~\ref{reflexive} in an early version of this paper. The
first author would like to thank members of the Department of
Mathematics and Statistics at Lancaster University and the
Department of Mathematics at the University of Iowa for kind
hospitality during the preparation of this article.

%%%%%%%%%%%%%%%%%%%%%%% REFERENCES
%%%%%%%%%%%%%%%%%%%%%%%%%%%%

%\begin{tabbing}
%{\it E-mail address}:xx\= \kill
%\noindent {\footnotesize\it Addresses}:
%\>{\footnotesize\sc Department of Mathematics}\\
%\>{\footnotesize\sc University of Iowa}\\
%\>{\footnotesize\sc Iowa City, IA\quad 52242}\\
%\>{\footnotesize\sc USA}\\
%\\
%\>{\footnotesize\sc Department of Mathematics and Statistics}\\
%\>{\footnotesize\sc Lancaster University}\\
%\>{\footnotesize\sc Lancaster, England}\\
%\>{\footnotesize\sc LA1 4YW}\\
%\\
%{\footnotesize\it E-mail addresses}:
%\>{\footnotesize\sf dkribs@math.uiowa.edu}\\
%\>{\footnotesize\sf s.power@lancaster.ac.uk}

%\end{tabbing}

\end{document}